\documentclass[a4paper,11pt]{amsart}
\usepackage{amssymb}
\usepackage{amsmath}
\usepackage{latexsym}
\usepackage{eepic}
\usepackage{epsfig}
\usepackage{graphicx}
\usepackage{amscd}
\usepackage{eufrak}
\usepackage{mathrsfs}
\usepackage[english]{babel}

\DeclareMathAlphabet{\mathpzc}{OT1}{pzc}{m}{it}
\newtheorem{theorem}[equation]{Theorem}
\newtheorem*{theorem*}{Theorem}
\newtheorem{theorem-definition}[equation]{Theorem-Definition}
\newtheorem{lemma-definition}[equation]{Lemma-Definition}
\newtheorem{definition-prop}[equation]{Proposition-Definition}

\newtheorem{prop}[equation]{Proposition}
\newtheorem*{prop*}{Proposition}
\newtheorem{lemma}[equation]{Lemma}
\newtheorem{cor}[equation]{Corollary}
\newtheorem{definition}[equation]{Definition}
\newtheorem*{definition*}{Definition}

\theoremstyle{definition}

\newtheorem{remark}[equation]{Remark}

\newcommand{\N}{\ensuremath{\mathbb{N}}}
\newcommand{\Z}{\ensuremath{\mathbb{Z}}}
\newcommand{\Q}{\ensuremath{\mathbb{Q}}}

\newcommand{\R}{\ensuremath{\mathbb{R}}}
\newcommand{\C}{\ensuremath{\mathbb{C}}}

\newcommand{\cX}{\ensuremath{\mathscr{X}}}
\newcommand{\mX}{\ensuremath{\mathfrak{X}}}
\newcommand{\mY}{\ensuremath{\mathfrak{Y}}}

\newcommand{\cA}{\ensuremath{\mathscr{A}}}

\newcommand{\cU}{\ensuremath{\mathscr{U}}}

\newcommand{\cC}{\ensuremath{\mathscr{C}}}
\newcommand{\cD}{\ensuremath{\mathscr{D}}}

\renewcommand{\R}{\ensuremath{\mathbb{R}}}
\renewcommand{\C}{\ensuremath{\mathbb{C}}}

\renewcommand{\cA}{\ensuremath{\mathscr{A}}}

\renewcommand{\cU}{\ensuremath{\mathscr{U}}}

\newcommand{\Spec}{\ensuremath{\mathrm{Spec}\,}}
\newcommand{\Spf}{\ensuremath{\mathrm{Spf}\,}}

\newcommand{\Lie}{\mathrm{Lie}}

\newcommand{\red}{\mathrm{red}}

\newcommand{\lcm}{\mathrm{lcm}}
\newcommand{\Jac}{\mathrm{Jac}}

\newcommand{\tame}{\mathrm{tame}}

\newcommand{\Gal}{\mathrm{Gal}}

\newcommand{\coker}{\mathrm{coker}}

\newcommand{\gp}{\mathrm{gp}}

\newcommand{\clog}{\omega^{\mathrm{log}}}
\newcommand{\dlog}{\mathrm{dlog}}
\newcommand{\Sch}{\mathrm{Sch}}
\newcommand{\sat}{\mathrm{sat}}

\numberwithin{equation}{subsection}

\newcommand{\sss}{\vspace{5pt} \subsubsection*{ }\refstepcounter{equation}{{\bfseries(\theequation)}\ }}

\begin{document}
\title[A logarithmic interpretation of Edixhoven's jumps]{A logarithmic interpretation of Edixhoven's jumps for Jacobians}

\author[Dennis Eriksson]{Dennis Eriksson}
\address{Gothenburg University and Chalmers Institute for Technology \\
Department of Mathematics\\
412 96 Gothenburg \\
Sweden}
\email{dener@chalmers.se}

\author[Lars Halvard Halle]{Lars Halvard Halle}
\address{University of Stavanger\\
Department of Mathematics and Natural Sciences\\
4036 Stavanger\\
Norway}
\email{lars.h.halle@uis.no}

\author[Johannes Nicaise]{Johannes Nicaise}
\address{KU Leuven\\
Department of Mathematics\\
Celestijnenlaan 200B\\3001 Heverlee \\
Belgium}
\email{johannes.nicaise@wis.kuleuven.be}

\begin{abstract}
 Let $A$ be an abelian variety over a discretely valued field. Edixhoven has defined a filtration on the special fiber of the N\'eron model of $A$ that
 measures the behaviour of the N\'eron model under tame base change. We interpret the jumps in this filtration in terms of lattices of logarithmic differential forms
 in the case where $A$ is the Jacobian of a curve $C$, and we give a compact explicit formula for the jumps in terms of the combinatorial reduction data of $C$.
\end{abstract}

\maketitle

\section{Introduction}
Let $R$ be a henselian discrete valuation ring, with quotient field $K$ and algebraically closed residue field $k$. We denote by $p$ the characteristic exponent of $k$.
 Let $A$ be an abelian variety over $K$, and denote by $\cA$ its N\'eron model; loosely speaking, this is the minimal smooth model of $A$ over $R$. If $K'$ is a finite separable extension of $K$ with valuation ring $R'$, then
  from the universal property that defines the N\'eron model we get a canonical morphism of $R'$-group schemes
  $h:\cA\times_R R' \to \cA'$ where $\cA'$ denotes the N\'eron model of $A'=A\times_K K'$. The morphism $h$ is not an isomorphism unless
   $K=K'$ or $A$ has good reduction (i.e., $\cA$ is an abelian $R$-scheme). One would like to understand the properties of the morphism $h$, and thus
    how the N\'eron model of $A$ behaves under base change.

 In \cite{edix}, Edixhoven constructed a filtration on the special fiber
 $\cA_k$  of $\cA$ by closed subgroups that measures the behaviour of the N\'eron model under {\em tame} extensions of the base field $K$, that is, finite extensions of degree prime to the characteristic exponent $p$ of $k$. This filtration is indexed by rational numbers in the semi-open interval $[0,1[$. The numbers in $[0,1[$ where this filtration jumps form an interesting set of invariants of $A$; we will simply refer to these numbers as the jumps of $A$. Note that, by definition, these jumps are {\em real} numbers. The question whether
 they are always rational is one of the main open problems about these invariants.

 Let $K^s$ be a separable closure of $K$. By Grothendieck's Semi-Stable Reduction Theorem \cite[IX.3.6]{SGA7-1} there exists a smallest finite extension
  $L$ of $K$ in $K^s$ such that the abelian variety $A\times_K L$ has semi-stable reduction; this means that the identity component of the special fiber of the N\'eron model
   of $A\times_K L$ is an extension of an abelian $k$-variety by an algebraic $k$-torus. We say that $A$ is tamely ramified if $L$ is a tame extension of $K$, and wildly ramified
   otherwise. When $A$ is tamely ramified, Edixhoven explained how the jumps of $A$ can be computed from the action of the Galois group $\Gal(L/K)$ on the N\'eron model of $A\times_K L$.
    This description implies in particular that the jumps are rational and that the degree of $L$ over $K$ is the least common multiple of their denominators.

  In this paper we are mainly interested in the wildly ramified case. We will study the jumps of $A$ under the assumption that $A$ is the Jacobian of a
  $K$-curve $C$. In this case, it was already proven by the second author in \cite{Halle-neron} that the jumps are rational; see also \cite[\S5.3]{HaNi-mem}. In fact, he proved a much stronger result, namely that the jumps of $A$
   only depend on the combinatorial reduction data of $C$, and not on the characteristic of $k$. A result of Winters \cite{winters} guarantees that we can always find a curve
   $D$ over the field of complex Laurent series $\C((t))$ with the same reduction data, and thus the same jumps, as $C$. Since $D$ is automatically tame, it follows that the jumps
   of $C$ are rational. By Corollary 3.1.5 in Chapter 5 of \cite{HaNi-mem}, the least common multiple of their denominators is the so-called {\em stabilization index} of $C$, a combinatorial invariant whose definition we will recall in \eqref{sss-stabcurve}.

    This result is quite powerful, but the proof of the rationality of the jumps is somewhat indirect, since one uses the combinatorial nature of the jumps and Winters's result to reduce to the case where the residue field $k$
    has characteristic zero. The aim of the present paper is to give an interpretation of the jumps of $A$ in terms of lattices of logarithmic differential forms on $R$-models of $C$ (Theorems \ref{thm-logtame} and \ref{thm-logdescr}), and to deduce a direct proof of their rationality (Theorem \ref{thm-rat}). Very roughly, we show that the jumps are controlled by the so-called {\em saturation} of a log regular model of $C$. The saturation is essentially equivalent to a semi-stable model in the tamely ramified case, but not in general; it can be viewed as a kind of combinatorial (characteristic-free) approximation of a semi-stable model.
       We expect that this approach can be generalized to arbitrary abelian $K$-varieties, provided that one finds a good notion of logarithmic N\'eron model of $A$. To our knowledge, such a construction has not yet appeared in the literature.
  Using our logarithmic interpretation of the jumps, we then establish an explicit
   and compact formula for the jumps of $A$ in terms of the combinatorial reduction data of $C$ (Theorem \ref{theorem-formula}).
   Such a formula was not known before.  It allows to prove directly that the stabilization index $e(C)$ is the least common multiple of the denominators of the jumps (Corollary \ref{cor-lcm}), without relying on Winters's result.
  We also deduce
    some interesting new properties of the jumps that make it easy to compute them in concrete examples (Propositions \ref{prop-lowbound} and  \ref{prop-princ}).
       Our methods give a conceptual explanation of the role that is played by the stabilization index of $C$: it is the smallest possible saturation index of a log regular model of $C$ (see Corollary \ref{cor-stab}). Our formula and the properties we deduce from it are somewhat reminiscent of jumping numbers of multiplier ideals of divisors on surfaces, but we do not know if there are any direct connections.

\subsection*{Acknowledgements}  A part of the research for this paper was done while JN
 was a member of the program {\em Model Theory, Arithmetic Geometry and Number Theory} at MSRI, Berkeley; he would like to thank the institute for its hospitality. JN was partially supported by the Fund for Scientific Research - Flanders (G.0415.10) and ERC Starting Grant MOTZETA. 

\subsection*{Notation}
 Throughout the paper, we will let $R$ denote a henselian discrete
valuation ring, with maximal ideal $\frak{m}$, quotient field $K$ and residue field $k$.
 We will assume that $ k $ is algebraically closed and we denote its characteristic exponent by $ p \geq 1
 $. We set $S=\Spec R$. Given a finite extension $K'$ of $K$, we
 will denote by $R'$ the integral closure of $R$ in $K'$ (which is
 again a henselian discrete valuation ring) and we set $S'=\Spec
 R'$.

 We fix a separable closure $ K^{s} $ of $K$, and we denote by
$K^t$ the tame closure of $K$ in $K^s$. The integral closures of
$R$ in $K^t$ and $K^s$ will be denoted by $R^t$ and $R^s$,
respectively.

 For every ring $A$ we denote by $(\Sch/A)$ the category of $A$-schemes. For every $A$-algebra $B$ we denote by
 $$(\cdot)_B:(\Sch/A)\to (\Sch/B):\cX\mapsto \cX_B = \cX \times_A B $$ the
 base change functor.

 Throughout the paper, $C$ will always denote a smooth, proper and geometrically connected
$K$-curve of genus $g > 0$. We also assume that $C$ has a zero divisor of degree one.
The Jacobian variety of $C$ is denoted $ \Jac(C) $.

 If $P$ is a monoid, then we write $P^{\gp}$ for its groupification and $P^{\mathrm{sat}}$ for its saturation.

\section{A few reminders on Edixhoven's filtration and N\'eron models of Jacobians}
\subsection{Chai's base change conductor}
\sss Let $A$ be an abelian $K$-variety of dimension $g$ and let
$\mathscr{A}$ denote its N\'eron model over $R$. Let moreover
$K'/K$ be a finite separable field extension and denote by $R'$
the integral closure of $R$ in $K'$ and by $\frak{m}'$ its maximal
ideal. We denote by $\mathscr{A}'$ the N\'eron model of $ A
\times_K K' $ over $R'$. Since $\mathscr{A} \times_R R'$ is smooth
and $\mathscr{A}'$ is a N\'eron model, there exists a unique
morphism
$$ h : \mathscr{A} \times_R R' \to \mathscr{A}' $$
extending the canonical isomorphism between the generic fibers. We
shall refer to this morphism as the \emph{base change morphism}.

\sss On the level of Lie algebras, the base change morphism
induces an injective homomorphism
$$\mathrm{Lie}(h) : \mathrm{Lie}(\mathscr{A}) \otimes_R R' \to \mathrm{Lie}(\mathscr{A}') $$
of free $R'$-modules of rank $g = \mathrm{dim}(A)$.

\begin{definition} \label{Def:basechange}\item
\begin{enumerate}
\item The tuple of $K'$-elementary divisors of $A$ is the unique
non-decreasing tuple $$(c_1(A,K'),\ldots,c_g(A,K'))$$ in $\N^g$
such that
$$\coker(\Lie(h))\cong \bigoplus_{i=1}^g \left(R'/(\frak{m}')^{c_i(A,K')}\right).$$
\item The tuple $(c_1(A),\ldots,c_g(A))$ of elementary divisors of
$A$ is defined by
$$c_i(A)=\frac{1}{[K':K]}c_i(A,K')$$ where $K'$ is any finite separable
extension of $K$ such that $A\times_K K'$ has semi-stable
reduction. The base change conductor $c(A)$ of $A$ is defined by
$$c(A)=\sum_{i=1}^g
c_i(A)=\frac{1}{[K':K]}\mathrm{length}_{R'}\coker(\Lie(h)).$$
\end{enumerate}
\end{definition}
\sss It follows from \cite[IX.3.3]{SGA7-1} that the definition of $c_i(A)$ and
$c(A)$ is independent of choice of the extension $K'/K$ over which
$A$ has semi-stable reduction. The base change conductor $c(A)$
 and the elementary divisors $c_i(A)$ were defined by Chai and Yu for algebraic tori in \cite{chai-yu}.
Chai generalized
 this definition to semi-abelian varieties in \cite{chai}. The base change conductor measures the defect of semi-stable
  reduction of $A$; in particular, it vanishes if and only if $A$
has semi-stable reduction.


\subsection{Edixhoven's filtration and the tame base change conductor}
\sss In \cite{edix}, Edixhoven defined a filtration
$F^{\bullet}\cA_k$  on the special fiber $\cA_k$ of the N\'eron
model $\cA$ of $A$, by closed subgroups $F^i\cA_k$ indexed by
rational numbers  $i\in [0,1[$. This filtration measures the
behaviour of the N\'eron model under tame finite extensions of
$K$. One can define the jumps in this filtration by looking at the
indices where $F^i\cA_k$ changes. These jumps are real numbers
 in $[0,1[$, and each of them has a natural multiplicity; see \cite[5.1.3.6]{HaNi-mem}. The
 number of jumps (counted with multiplicities) is equal to the
 dimension $g$ of the abelian variety $A$.

\sss Edixhoven also proposed an alternative way to compute the
 jumps and their multiplicities, which we can reformulate as follows.
 Let $K_0\subset K_1\subset \ldots$ be a tower of finite
 extensions of $K$ in $K^t$ that is cofinal in the set of all
 finite extensions of $K$ in $K^t$, ordered by inclusion. Then for
 every index $i$ in $\{1,\ldots,g\}$, the sequence
 $$(c_i(A,K_n)/[K_n:K])_{n\in \N}$$ is a non-decreasing sequence of
 rational numbers; one can show that the elements of this sequence are strictly bounded
 by $1$. We denote by $j_i(A)$ the limit of this sequence.
  By \cite[5.1.3.7]{HaNi-mem}, the non-decreasing tuple of real numbers
  $$(j_1(A),\ldots,j_g(A))$$ is precisely the tuple of jumps of
  $A$, counted with multiplicities.

 \sss In \cite[5.1.3.6]{HaNi-mem}, the last two authors defined the
{\em tame base change conductor} $c_{\tame}(A)$ of $A$ as the sum
of the jumps:
$$c_{\mathrm{tame}}(A)=\sum_{i=1}^g j_i(A).$$
 If
$A$ is tamely ramified, then $j_i(A)=c_i(A)$ for all $i$ and
$c_{\tame}(A)=c(A)$
 by \cite[4.18]{HaNi}, but these invariants differ in general.

 \sss Whereas the elementary divisors $c_i(A)$ are rational numbers by definition, it is an open problem whether the jumps
 $j_i(A)$ are always rational numbers (this question was already raised by Edixhoven in \cite[5.4.5]{edix}). It is not even known whether their sum $c_{\mathrm{tame}}(A)$ is always
 rational. Assuming that all the jumps are rational, we define the
 stabilization index $e(A)$ of $A$ as the smallest positive
 integer $e$ such that $e\cdot j_i(A)$ is rational for every $i$.
  If $A$ is tamely ramified, it follows from
 the equalities $c_i(A)=j_i(A)$ that the jumps are rational numbers, and one
 can show that $e(A)$ is equal to the degree of the
 minimal extension $L$ of $K$ in $K^s$ such that $A\times_K L$
 has semi-stable reduction \cite[5.1.3.12]{HaNi-mem}.

\sss \label{sss-stabdef} The stabilization index of $A$ (whose definition depends on
the rationality of the jumps) seems to capture important
information about the behaviour of the N\'eron model of $A$ under
 tame base change: the idea is that the N\'eron model should
 change ``as little as possible'' if the degree of the base
 extension is prime to $e(A)$. This principle is made precise in
 section 1 of Part 4 in \cite{HaNi-mem}; it is the key to
  understanding the so-called motivic zeta function of $A$.

\subsection{Regular models of curves}
\sss Let $C$ be a smooth, proper and geometrically connected
$K$-curve of genus $g > 0$. We assume that $C$ has a zero divisor of degree one.
 An $R$-model of $C$ is a flat and proper $R$-scheme $\cC$,
 endowed with an isomorphism of $K$-schemes
$$ \cC \times_R K \to C. $$ Morphisms of models are defined in the obvious way.
  An $ncd$-model (resp.~$sncd$-model) of $C$ is a regular $R$-model $\cC$
 such that the special fiber $\cC_k$ is a divisor with normal
 crossings (resp.~strict normal crossings, i.e. the reduced irreducible components are smooth). Since we assume that
 the genus of $C$ is at least one, the curve $C$ has a unique
 minimal regular model, a unique minimal $ncd$-model and a unique
 minimal $sncd$-model; we refer to section 2.2 of \cite{ni-saito}
 for a brief summary of the literature with detailed references. We say that $\cC$ has semi-stable reduction
 if the special fiber of the minimal $ncd$-model of $C$ is reduced. By the Semi-Stable Reduction Theorem for curves \cite{deligne-mumford}, there exists a smallest finite extension $L$ of $K$
  in $K^s$ such that $C\times_K L$ has semi-stable reduction. Moreover, since we assume that $C$ has a zero divisor of degree one, $C$ has semi-stable reduction if
   and only if  its Jacobian variety has semi-stable reduction (this does not hold for genus one curves without a rational point, but such a curve
    never has a zero divisor of degree one).

\sss \label{sss-cohflat} If $\cC$ is an $sncd$-model of $C$ with special fiber $\cC_k=\sum_{i=1}^{r}N_iE_i$, then the {\em combinatorial reduction data} of $\cC$ consist of the dual graph of $\cC_k$, where we label the vertex corresponding to $E_i$ with the multiplicity $N_i$ and the genus $g(E_i)$ of $E_i$. We define the combinatorial reduction data of $C$ to be those of the minimal $sncd$-model of $C$.
 The condition that $C$ has a zero divisor of degree one is equivalent to the condition that the greatest common divisor
 of the multiplicities $N_i$ is equal to one, by \cite[7.1.6]{raynaud}. This implies that the structural morphism $\cC\to S$ is cohomologically flat \cite[8.2.1]{raynaud}, so that $H^1(\cC,\mathcal{O}_{\cC})$ is a free $R$-module.

\sss \label{sss-stabcurve} If $\cC$ is an $sncd$-model of $C$, then we say that an irreducible component $E$ of the special fiber
 $\cC_k$ is {\em principal} if the genus of $E$ is at least one or $E$ is a rational curve meeting the other components of $\cC_k$ in at least three points.
  The {\em stabilization index} of $C$ is defined as the least common multiple of the multiplicities of the principal components in the special fiber of the
   minimal $sncd$-model of $C$; see Definition 2.2.2 in Chapter 3 of \cite{HaNi-mem}. Let $L$ be the minimal finite extension
 of $K$ in $K^s$ such that $C\times_K L$ has semi-stable reduction. If this extension is tame, then its degree is equal to $e(C)$ by \cite[3.4.4]{ni-saito}, but this is false in general. As we have already mentioned in the introduction, it was proven by the last two authors   that the stabilization index of $C$ is equal to that of its Jacobian in the sense of \eqref{sss-stabdef} (see Corollary 3.1.5 in Chapter 5 of \cite{HaNi-mem}). We will give a more direct and conceptual argument in Corollary \ref{cor-lcm}.

\subsection{N\'eron models of Jacobians}
\sss Let $\cC$ be a regular model of $C$, let $ \mathscr{J} $ be
the N\'eron model of $J=\Jac(C)$ and let $ \mathscr{J}^0 $ be its
identity component. Recall that we assume that the curve $C$ has
 a zero divisor of degree one.  By a fundamental theorem of Raynaud
\cite[9.5.4]{neron}, there is a natural isomorphism
$$ \mathrm{Pic}^0_{\cC/R} \cong \mathscr{J}^0.$$
 Via this description of $\mathscr{J}^0$, it is
possible to reduce many computations concerning N\'eron models to
computations on regular models of curves, something which is often
very useful. We will see that this is true, in particular, for
 the jumps in Edixhoven's filtration. To this end, we will need the
following well known interpretation of the invariant differential forms on
$\mathscr{J}$.

\sss  Let $ e_{\mathscr{J}} : S \to \mathscr{J} $ be the unit
section of the N\'eron model of $J$. Then we define $\Omega(J)$ to
be the module of translation-invariant relative differential forms
on $\mathscr{J}$ (see \cite[\S4.2]{neron}); thus
$$ \Omega(J)= e_{\mathscr{J}}^* \Omega^1_{\mathscr{J}/S}. $$
 This is an $R$-lattice in the $g$-dimensional $K$-vector space
 $$H^0(J,\omega_{J/K})\cong H^0(C,\omega_{C/K}).$$

\begin{prop}\label{prop-compjac}
Let $\cC$ be a regular $R$-model of $C$, and denote by
$\omega_{\cC/R}$ the relative canonical sheaf. Then one has
$$ \Omega(J) = H^0(\cC,\omega_{\cC/R}) $$
as lattices in $H^0(C,\omega_{C/K})$.
\end{prop}
\begin{proof}
By \cite[8.4.1]{neron}, the canonical isomorphism
$$\mathrm{Lie}(J)\to H^1(C,\mathcal{O}_{C})$$ can be extended to
an isomorphism of $R$-modules
$$\mathrm{Lie}(\mathscr{J})\to H^1(\cC,\mathcal{O}_{\cC}).$$ Since
$\omega_{\cC/R}$ is a dualizing sheaf for the structural morphism
$\cC\to S$, Grothendieck duality yields an isomorphism
$$\mathrm{Lie}(\mathscr{J})^{\vee}\cong \Omega(J)\to H^0(\cC,\mathcal{\omega}_{\cC/R}).$$
\end{proof}

\section{Log schemes and differential forms}
In this section we will prove some basic results on sheaves of differentials on logarithmic schemes.
 The standard introduction to logarithmic geometry is \cite{kato-log}.
  Our log structures are defined with respect to the \'etale topology, and our definition of log regularity is the one from \cite[2.2]{niziol}; see
  \cite[2.3]{niziol} for a comparison with Kato's definition for the Zariski topology in \cite{kato-toric}.
\subsection{Base change of $fs$ log schemes}
\sss We will denote by $S^+$ the scheme $S=\Spec(R)$ with its
standard log structure, that is, the log structure induced by
 the morphism of monoids $R\setminus \{0\}\to R$. If $\cX$ is a flat
 $R$-scheme, we will denote by $\cX^+$ the log scheme we
 obtain by endowing $\cX$ with the divisorial log structure
 associated with $\cX_k$. Similar notation will be used for
 $S'=\Spec(R')$ and schemes over $R'$, if $K'$ is a finite
 extension of $K$ and $R'$ is its valuation ring. The fiber
 product in the category of fine and saturated ($fs$) log schemes will be denoted
 by $\times^{fs}$. Beware that it does not commute with the
 forgetful functor to the category of schemes.

\sss  Let $\cC$ be an $R$-model of $C$ such that $\cC^+$ is log
regular. Let $K'$ be a finite extension of $K$.
 Set $$\mathscr{D}^+=\cC^+\times^{fs}_{S^+}(S')^+.$$  Let
 $\mathscr{D}$ be the underlying scheme of $\mathscr{D}^+$. Then the log structure on $\mathscr{D}^+$ is the divisorial log structure induced by
  $\mathscr{D}_k$, so that our notation is consistent. The scheme  $\mathscr{D}$ is
 canonically equipped with a finite morphism
 $$\mathscr{D}\to \mathscr{C}\times_R R'$$ which is an isomorphism
 on the generic fibers since there the log structure is trivial.
\sss \label{log-smooth} If $\cC^+$ is log smooth over $S^+$, then  $\mathscr{D}^+$ is
log smooth over  $(S')^+$, because log smoothness is preserved by
 base change in the category of $fs$ log schemes. Likewise, if $K'$
 is a tame extension of $K$, then $(S')^+$ is log \'etale over
 $S^+$ and $\cD^+$ is log \'etale over $\cC^+$. In both cases,
 $\cD^+$ is log smooth over a log regular scheme, and thus itself log regular \cite[8.2]{kato-toric}, which implies that the
 underlying scheme $\mathscr{D}$ is normal \cite[4.1]{kato-toric}. Therefore,
 $$\mathscr{D}\to \mathscr{C}\times_R R'$$ must be a normalization
 map.


\subsection{Differential forms on log regular schemes}\label{subsec-log}

\sss Let $\cX$ be a flat $R$-scheme of finite type such that the associated log scheme $\cX^+$ is log regular and let $x$ be a closed point of $\cX_k$.
 The aim of this section is to construct a free resolution of the module of germs of log differential forms $\Omega^1_{\cX^+/S^+,x}$.
   This result will be important for the study of logarithmic canonical sheaves in the following sections. To simplify the notations, we introduce the following
   abbreviations: we write $\mathcal{O}$ for the \'etale-local ring of $\cX$ at $x$ (the henselization of $\mathcal{O}_{\cX,x}$), $\mathcal{M}=\mathcal{O}\cap (\mathcal{O}\otimes_R K)^{\times}$ for the
    monoid  $\mathcal{M}_{\cX^+,x}$,
       $\overline{\mathcal{M}}=\mathcal{M}/\mathcal{O}^{\times}$ for the characteristic monoid and $\Omega$ for the $\mathcal{O}$-module
  $\Omega^1_{\cX^+/S^+,x}$ (where the stalk at $x$ is taken in the \'etale topology).
   Recall that the monoid $\overline{\mathcal{M}}$ is toric (that is, sharp, fine and saturated).
  We fix a section $\overline{\mathcal{M}}\to \mathcal{M}$ for the projection morphism $\mathcal{M}\to \overline{\mathcal{M}}$, so that we can view
 $\overline{\mathcal{M}}$ as a submonoid of the multiplicative monoid $(\mathcal{O},\cdot)$.  Denote by $I$ the ideal of $\mathcal{O}$ generated by $\overline{\mathcal{M}}\setminus \{1\}$.
   By definition of log regularity, the local ring $\mathcal{O}/I$ is regular, and its dimension $r$ is equal to the dimension of $\mathcal{O}$ minus the rank of
    the free abelian group $\overline{\mathcal{M}}^{\gp}$.  We choose elements $t_1,\ldots,t_r$ in $\mathcal{O}$ whose reductions modulo $I$ form a regular system of local parameters in $\mathcal{O}/I$, and we denote by $M$ the $\mathcal{O}$-module $$M=\mathcal{O}\otimes_{\Z}(\overline{\mathcal{M}}^{\gp}\oplus \Z^r).$$

\begin{prop}\label{prop:freeres}
Denote by
 $$\gamma:M\to \Omega$$
  the unique morphism of $\mathcal{O}$-modules that sends $a\otimes (m\oplus n)$ to $$a(\dlog(m)+\sum_{i=1}^r n_idt_i)$$ for all $a$ in $\mathcal{O}$, all $m$ in
  $\overline{\mathcal{M}}$ and all $n$ in $\Z^r$. Then $\gamma$ is surjective, and its kernel is a free $\mathcal{O}$-module of rank one.
\end{prop}
\begin{proof}
By faithful flatness of the completion morphism $\mathcal{O}\to \widehat{\mathcal{O}}$, it is enough to prove the statement after base change to
$\widehat{\mathcal{O}}$. We denote by $\widehat{R}$ the completion of $R$ and by $\mX$ the formal $\widehat{R}$-scheme $\Spf \widehat{\mathcal{O}}$, and we define a log structure on $\mX$ by means of the chart
$\overline{\mathcal{M}}\to \widehat{\mathcal{O}}$. The resulting log formal scheme will be denoted by $\mX^+$, and we write $\mathfrak{S}^+$ for the formal scheme $\Spf \widehat{R}$ with its
 standard log structure. Then we can identify $\Omega\otimes_{\mathcal{O}}\widehat{\mathcal{O}}$ with the module of log differentials $\Omega^1_{\mX^+/\mathfrak{S}^+}$.
 Set $A=\widehat{R}[[\overline{\mathcal{M}}]][[T_1,\ldots,T_r]]$ and denote by $\mY^+$ the log formal scheme $\Spf A$ with chart $\overline{\mathcal{M}}\to A$.
   By the local description of toric singularities in \cite[3.2]{kato-toric}, we know that we can view $\mX^+$ as a strict closed log formal subscheme of
   $\mY^+$ defined by a principal ideal $J$ such that $t_i$ is the restriction of $T_i$ to $\mX^+$ for every $i$ (in the notation of \cite[3.2(2)]{kato-toric} the ring $R$ is the ring of
    Witt vectors $W(k)$, but the proof can be adapted in an obvious way). Let $m_1,\ldots,m_s$ be a basis for $\overline{\mathcal{M}}^{\gp}$. Then  $\Omega^1_{\mY^+/\mathfrak{S}^+}$ is free with basis
   $dT_1,\ldots,dT_r,\dlog(m_1),\ldots,\dlog(m_s)$. Thus the base change of $\gamma$ to $\widehat{\mathcal{O}}$ fits into the fundamental short exact sequence of $\widehat{\mathcal{O}}$-modules
   $$0\to J/J^2 \to \Omega^1_{\mY^+/\mathfrak{S}^+}\otimes_{A}(A/J)=M\otimes_{\mathcal{O}}\widehat{\mathcal{O}}\to \Omega^1_{\mX^+/\mathfrak{S}^+}=\Omega\otimes_{\mathcal{O}}\widehat{\mathcal{O}}\to 0. $$
\end{proof}

\subsection{Logarithmic canonical sheaves}\label{subsec-logcan}

\sss Let $X$ be a Noetherian scheme and $\mathcal{F}$ a coherent
$\mathcal{O}_{X}$-module. Recall that the {\em reflexive hull} of
$\mathcal{F}$ is the double dual $\mathcal{F}^{\vee \vee}$ of
$\mathcal{F}$, and that $\mathcal{F}$ is called {\em reflexive} if
the natural morphism $\mathcal{F}\to \mathcal{F}^{\vee \vee}$ is
an isomorphism. We recall a few basic properties of
reflexive sheaves:
\begin{itemize}
\item If $X$ is regular, then every reflexive rank one sheaf on
$X$ is a line bundle \cite[1.9]{har80}. \item If $X$ is
normal, then every reflexive sheaf $\mathcal{F}$ on $X$ has the
$S_2$ property \cite[1.9]{har94}. This implies that, for
every closed subscheme $Z$ of $X$ of codimension at least two, the
restriction map $$\mathcal{F}(X)\to \mathcal{F}(X\setminus Z)$$ is
an isomorphism. \item Assume that $X$ is normal and let $Z$ be a closed subscheme of $X$ of
codimension at least two. Let $\mathcal{G}$ be a reflexive sheaf on
$U=X\setminus Z$. If we denote by $i$ the open immersion $U\to X$,
then $i_*\mathcal{G}$ is a reflexive sheaf on $X$, and it is the
unique extension of $\mathcal{G}$ to a reflexive sheaf on $X$.
\end{itemize}

\sss Let $\mathscr{C}$ be a normal $R$-model of $C$. Denote by
$\cU=\cC^{\mathrm{reg}}$ the open subscheme of regular points of
$\cC$ and by $i:\cU\to \cC$ the open immersion. We denote by
$\omega_{\cU/R}$ the canonical line bundle of the morphism $\cU\to
\Spec R$ and we define the canonical sheaf of the $R$-scheme $\cC$
by
$$\omega_{\cC/R}=i_*\omega_{\cU/R}.$$
 This is a reflexive rank one sheaf on $\cC$, whose restriction to $C$
is naturally isomorphic to the canonical line bundle
$\omega_{C/K}$. If the structural morphism $f:\cC\to \Spec R$ is
l.c.i., then $\omega_{\cC/R}$ is canonically isomorphic to the
relative canonical bundle of $f$, i.e., the determinant of
$\Omega^1_{\cC/R}$. We say that $\cC$ has canonical singularities if, for every morphism $g:\cC'\to \cC$ of $R$-models
 of $C$ such that $\cC'$ is regular, we have $g^*\omega_{\cC/R}\subset \omega_{\cC'/R}$ as subsheaves of
  $j_*\omega_{C/K}$, where $j$ denotes the open immersion $C\to \cC'$.

\sss Likewise, we define the logarithmic canonical sheaf on $\cC$
by
$$\clog_{\cC/R}=i_*(\det \Omega^1_{\cU^+/S^+}).$$  This is again a
 reflexive rank one sheaf on $\cC$ whose restriction to $C$
is naturally isomorphic to the canonical line bundle
$\omega_{C/K}$. Its relation to $\omega_{\cC/R}$ is explained in
the following proposition.

\begin{prop}\label{prop-compar}
Let $\cC$ be a normal model of $C$, and denote by $j$ the open
immersion $j:C\to \cC$. Then $$\clog_{\cC/R}
=\omega_{\cC/R}(\cC_{k,\red}-\cC_k)$$ as subsheaves of
$j_*\omega_{C/K}$. In particular, if $\cC_k$ is reduced, then
$\omega_{\cC/S}$ and $\clog_{\cC/R}$ coincide.
\end{prop}
\begin{proof}
 Since both sheaves are reflexive, it suffices to prove that they
 coincide on the complement of a finite set of closed points, so
 that we can replace $\cC$ by a regular open subscheme $\cU$ such
 that $\cU_k$ has strict normal crossings.
 Then the statement
    follows from \cite[5.3.4]{kato-saito} by taking
 determinants.
\end{proof}

\sss \label{sss-ramlog} The logarithmic canonical sheaf
$\clog_{\cC/R}$ behaves well under  $fs$ base change, in the following sense. Assume either that $\cC^+$ is log smooth over $S^+$, or that $\cC^+$ is log regular and $K'$
is a tame finite extension of $K$. If we set
$$\cD^+=\cC^+\times^{fs}_{S^+}(S')^+$$ then
 $\clog_{\cD/R'}$ is canonically isomorphic to the pullback of
 $\clog_{\cC/R}$ to $\cD$. This follows from the
 fact that $\cD$ is normal (cf. \ref{log-smooth}),
 $\Omega^1_{\cD^+/(S')^+}$ is the pullback of
 $\Omega^1_{\cC^+/S^+}$ and $\cD\to \cC$ is flat at every point of codimension $\leq
1$ (recall that taking determinants
commutes with flat base change).

\begin{prop}\label{prop-linebundle}
Let $\cC$ be an $R$-model of $C$ such that $\cC^+$ is log regular. Then the logarithmic canonical sheaf
 $\clog_{\cC/R}$ is the determinant line bundle of the perfect coherent sheaf $\Omega^1_{\cC^+/S^+}$. If $h:\cD\to \cC$ is a morphism of models of $C$ such that
 the morphism of log schemes $\cD^+\to \cC^+$ is log \'etale, and if we denote by $j$ the open immersion
 $\cD_K\to \cD$, then we have $\clog_{\cD/R}=h^*\clog_{\cC/R}$ as subsheaves of $j_*\omega_{C/K}$.
\end{prop}
\begin{proof} The entire statement is local with respect to the \'etale topology on $\cC$. Thus we may assume, by
 Proposition \ref{prop:freeres}, that
 there exists a resolution of $\Omega^1_{\cC^+/S^+}$ by free coherent sheaves of the form
 \begin{equation}\label{eq-res}
 0\to \mathcal{O}_{\cC}\to \mathcal{F}\to \Omega^1_{\cC^+/S^+}\to 0 .
 \end{equation}
 The determinant line bundle $\det(\Omega^1_{\cC^+/S^+})$ is equal to $\clog_{\cC/R}$ because these reflexive sheaves coincide on the regular locus of $\cC$ and $\cC$ is normal.
  Since $h$ is log \'etale, the pullback $h^*\Omega^1_{\cC^+/S^+}$ is isomorphic to $\Omega^1_{\cD^+/S^+}$. Applying the right exact functor
  $h^*$ to our resolution \eqref{eq-res} yields a right exact sequence
 $$ \mathcal{O}_{\cD}\to h^*\mathcal{F}\to \Omega^1_{\cD^+/S^+}\to 0.$$ This sequence is also exact on the left because
  this holds on the generic fiber $C$ of $\cD$ and $\cD$ is flat over $R$. Thus
  $$\clog_{\cD/R}=\det(\Omega^1_{\cD^+/S^+})=h^*\clog_{\cC/R}. $$
\end{proof}

\section{A logarithmic interpretation of the jumps in Edixhoven's filtration}

\subsection{Comparing lattices over discrete valuation rings}
\sss Let $V$ be a vector space of dimension $g$ over $K$. For
every pair of $R$-lattices $L_0\subset L_1$ in $V$, we define the
tuple of {\em elementary divisors} of this pair as the unique
non-decreasing tuple
$$(c_1(L_1/L_0),\ldots,c_g(L_1/L_0))$$ in $\N^g$ such
that
$$L_1/L_0\cong \oplus_{i=1}^g R/\frak{m}^{c_i(L_1/L_0)}.$$ The {\em
conductor} $c(L_1/L_0)$ of the pair of lattices is defined by
$$c(L_1/L_0)=\sum_{i=1}^g
c_i(L_1/L_0)=\mathrm{length}_{R}(L_1/L_0).$$

\sss  Denote by $R^s$ the valuation ring of $K^s$. Let $L_0\subset L_1$ be a pair of $R^s$-lattices in
$V\otimes_K K^s$. Then we can choose a finite extension $K'$ of
$K$ in $K^s$ and $R'$-lattices $L'_0\subset L'_1$ in $V\otimes_K
K'$ such that $L_i=L'_i\otimes_{R'}R^s$ for $i=0,1$. We define the
 elementary divisors $c_1(L_1/L_0),\ldots,c_g(L_1/L_0)$ and the conductor
$c(L_1/L_0)$ of the pair $(L_0,L_1)$ by
\begin{eqnarray*}
c_i(L_1/L_0)&=&\frac{1}{[K':K]}c_i(L'_1/L'_0),
\\c(L_1/L_0)&=&\frac{1}{[K':K]}c(L'_1/L'_0).
\end{eqnarray*}
It is straightforward to check that these definitions do not
depend on the choice of $K'$.
 We will make use of the following elementary property.

\begin{prop}\label{prop-elementary}
 Let $L_0\subset L_1\subset L_2$ be lattices in $V\otimes_K K^s$, and let $a$ be a non-zero element of $R^s$ such that the $R$-module $L_1/L_0$ is killed by
$a$.
 If we denote by $N$ the valuation of $a$ in $R^s$ (with respect to the unique valuation on $R^s$ that extends the normalized discrete valuation on $R$),
  then
$$c_i(L_2/L_1)\leq c_i(L_2/L_0)\leq c_i(L_2/L_1)+N$$ for every $i$ in
$\{1,\ldots,g\}$.
\end{prop}
\begin{proof}
 We can assume that the lattices $L_0$, $L_1$ and $L_2$ are defined over $R$, and that $a$ is an element of $R$. The proof is based on the following elementary
observation: for every positive integer $M$, the number of
elementary divisors of the pair $(L_0,L_2)$ greater than or equal
to $M$ is equal to the dimension of the $k$-vector space
 $$V_{0,2}^M=(\frak{m}^{M-1}L_2/L_0)\otimes_R k$$ and the analogous statement holds for
  $(L_1,L_2)$. The projection $L_2/L_0\to L_2/L_1$ gives rise to a surjection
  $V_{0,2}^M\to V_{1,2}^M$ for every $M$, which means that $c_i(L_2/L_1)\leq c_i(L_2/L_0)$ for all $i$.
    On the other hand, multiplication with $a$ defines a morphism of $R$-modules $L_2/L_1\to L_2/L_0$ which induces
    a surjection $V_{1,2}^M\to V_{0,2}^{M+N}$ for every $M$, so that we also have $c_i(L_2/L_0)\leq c_i(L_2/L_1)+N$ for all $i$.
\end{proof}

\subsection{Lattices of differential forms}
\begin{prop}\label{prop-indep}\item
Let $\cC$ be a normal model of $C$ over $R$.
\begin{enumerate}
\item \label{it:logindep} If $\cC^+$ is log regular, then the
$R$-lattice
$$H^0(\cC,\clog_{\cC/R})$$ in
$H^0(C,\omega_{C/K})$ only depends on $C$, and not on the choice
of $\cC$. \item \label{it:canindep} If we assume that $\cC$ has at
worst canonical singularities, then the $R$-lattice
$$H^0(\cC,\omega_{\cC/R})$$ in
$H^0(C,\omega_{C/K})$ only depends on $C$, and not on the choice
of $\cC$.

\end{enumerate}
\end{prop}
\begin{proof}
\eqref{it:logindep} Let $\cC_1$ be an $R$-model of $C$ such that $\cC_1^+$ is log regular. Let $f:\cC_2\to \cC_1$ be a morphism of $R$-models of $C$ that is obtained by blowing up $\cC_1$ at a closed point $x$ of
 its special fiber. Since $\cC_1$ has only rational singularities, the scheme $\cC_2$ is normal \cite[1.5]{lipman}.
  It suffices to show that  $\cC_2^+$ is log regular and
$$H^0(\cC_1,\clog_{\cC_1/R})=H^0(\cC_2,\clog_{\cC_2/R}),$$ since any
 pair of log regular $R$-models can be connected by a chain of such point blow-ups.

 It is clear that
$$H^0(\cC_2,\clog_{\cC_2/R})\subset H^0(\cC_1,\clog_{\cC_1/R})$$
 because $f$ is an isomorphism over $\cC_1\setminus \{x\}$ and
 $\clog_{\cC_1/R}$ is reflexive. Thus it is enough to prove that
\begin{equation}\label{eq-incl}f^*\clog_{\cC_1/R}\subset
\clog_{\cC_2/R}\end{equation} as subsheaves of the pushforward of $\omega_{C/K}$ to $\cC_2$.

 First, assume that $x$ is a regular point of the reduced special
 fiber $(\cC_{1,k})_{\red}$ of $\cC_1$. Then $\cC_1$ is also
 regular at $x$ by \cite[5.2]{niziol}. Thus $\cC^+_2$ is log regular and \eqref{eq-incl} is a straightforward consequence
 of the analogous inclusion for the canonical sheaves
 $\omega_{\cC_i/R}$, together with Proposition \ref{prop-compar}.

 Now assume that $x$ is a singular point of $(\cC_{1,k})_{\red}$.
  Then the morphism $\cC_2^+\to \cC_1^+$ is a log blow-up by \cite[4.3]{niziol},
 and hence log-\'etale, so that $\cC^+_2$ is log regular and \eqref{eq-incl} follows from Proposition \ref{prop-linebundle}.

 Property \eqref{it:canindep} can be proven in a similar way,
 replacing $\cC_1$ by a normal model of $C$ with at worst
 canonical singularities and $f:\cC_2\to \cC_1$ by a resolution of
 singularities of $\cC_1$; then the inclusion
 $$f^*\omega_{\cC_1/R}\subset \omega_{\cC_2/R}$$ follows from the
 definition of a canonical singularity.
\end{proof}

\begin{definition}\label{def-latt}
Let $\cC$ be a normal $R$-model of $C$. If $\cC^+$ is log regular, then we call the
$R$-lattice
$$H^0(\cC,\clog_{\cC/R})\subset H^0(C,\omega_{C/K})$$ the logarithmic lattice associated with
$C$, and we denote it by $\Omega_{\log}(C)$. If $\cC$ has at worst
canonical singularities, then we call the $R$-lattice
$$H^0(\cC,\omega_{\cC/R})\subset H^0(C,\omega_{C/K})$$
the canonical lattice associated with $C$, and denote it by
$\Omega_{\mathrm{can}}(C)$.
\end{definition}

\sss \label{sss-canlog} We can always find an $R$-model $\cC$ of $C$ such that $\cC^+$ is log regular (for instance, an $ncd$-model).
  By Proposition \ref{prop-indep}, Definition \ref{def-latt} does not depend on the choice of $\cC$.
 Note that
$$\Omega_{\mathrm{log}}(C)\subset\Omega_{\mathrm{can}}(C)$$ by
Proposition \ref{prop-compar}, and that they are equal when
 $C$ has
 semi-stable reduction.

\begin{prop}\label{prop-lattincl}
 Let $K'$ be a finite separable extension of $K$ and denote by
 $R'$ the integral closure of $R$ in $K'$. Set $C'=C\times_K K'$.
 \begin{enumerate}
\item \label{it:lattinc1.1} We have $$
\Omega_{\mathrm{can}}(C)\otimes_R
 R'\supset \Omega_{\mathrm{can}}(C')
$$
  as lattices in $H^0(C',\omega_{C'/K'})$, with equality if $C$ has semi-stable reduction.
  \item \label{it:lattinc1.3} Assume either that $C$ has an $R$-model $\cC$ such that $\cC^+$ is log smooth over $S^+$, or that $K'$ is a tame extension of $K$. Then we have $$
 \Omega_{\mathrm{log}}(C)\otimes_R
 R'\subset \Omega_{\mathrm{log}}(C')
$$
  as lattices in $H^0(C',\omega_{C'/K'})$.
 \end{enumerate}
\end{prop}
\begin{proof}
\eqref{it:lattinc1.1} The inclusion follows from from Proposition \ref{prop-compjac}. If $\cC$ is a semi-stable $R$-model of $C$,
 then
 $\cC\times_{R}R'$ is a normal $R'$-model of $C'$ with canonical
 singularities (rational double points of type $A_n$) so that
 $$\Omega_{\mathrm{can}}(C')=\Omega_{\mathrm{can}}(C)\otimes_{R} R'.$$

\eqref{it:lattinc1.3} This follows from \eqref{sss-ramlog}.
\end{proof}

\sss \label{sss-diag} In other words, the logarithmic lattice grows under tame extensions of $K$, and the canonical lattice shrinks under arbitrary extensions of $K$.
  It will be convenient to summarize all of the above inclusions in the following diagram; here $K'$ is a finite extension of $K$ in $K^t$, $K''$ is a
 finite extension of $K'$ in $K^s$, and $R'$ and $R''$ denote their respective valuation rings.
 $$\begin{CD}
 \Omega_{\mathrm{can}}(C)\otimes_R R^s &\ \supset\ & \Omega_{\mathrm{can}}(C\times_K K')\otimes_{R'} R^s &\ \supset\ & \Omega_{\mathrm{can}}(C\times_K K'')\otimes_{R''} R^s
 \\ \cup & & \cup & &
 \\  \Omega_{\mathrm{log}}(C)\otimes_R R^s &\ \subset\ & \Omega_{\mathrm{log}}(C\otimes_K K')\otimes_{R'} R^s & &
 \end{CD}$$

Now we can associate with the curve $C$ two further lattices.

 \sss The {\em semi-stable lattice} is defined by
$$ \Omega_{\mathrm{ss}}(C)=\Omega_{\mathrm{can}}(C\times_K K')\otimes_{R'} R^s\subset
H^0(C,\omega_{C/K})\otimes_K K^s,$$
 where $K'$ is any finite extension of $K$ in $K^s$ such that
 $C\times_K K'$ has semi-stable reduction. It follows from Proposition \ref{prop-lattincl} that this definition does not depend on the choice of
 $K'$, and
that $\Omega_{\mathrm{ss}}(C)$ is the intersection of the lattices
$$\Omega_{\mathrm{can}}(C\times_K K')\otimes_{R'} R^s$$ where $K'$
ranges over all the finite extensions of $K$ in $K^s$.

\sss The {\em saturated
 lattice}
 $\Omega_{\mathrm{sat}}(C)$
 is defined by $$
\Omega_{\mathrm{sat}}(C)= \bigcap_{K'}(\Omega_{\mathrm{can}}(C\times_K K')\otimes_{R'} R^s)
 $$  where $K'$ runs through the finite extensions of $K$ in $K^t$.
  We will prove in
 Theorem \ref{thm-logtame} that this is indeed a lattice in $H^0(C,\omega_{C/K})\otimes_K
 K^s$.
\if false

\begin{prop}\label{prop-inter}We have \begin{equation}\label{eq-inter}
\Omega_{\mathrm{sat}}(C)= \bigcap_{K'}(\Omega_{\mathrm{can}}(C\times_K K')\otimes_{R'} R^s)
 \end{equation}   where $K'$ runs through the finite extensions of $K$ in $K^t$, and $R'$ denotes the integral closure of
 $R$ in $K'$.
\end{prop}
\begin{proof}
The inclusion $\subset$ is a consequence of \eqref{sss-canlog} and Proposition \ref{prop-lattincl}. The converse inclusion $\supset$ follows from Lemma \ref{lemm-bounded}, applied
 to the base changes of $C$ to all tame finite extensions $K'$ of $K$, which tells us that we can approximate each element in the right hand side of  \eqref{eq-inter}
  with arbitrary precision by an element of  $\Omega_{\mathrm{log}}(C\times_K K')\otimes_{R'} R^s$ for some finite extension $K'$ of $K$ in $K^t$
    (note that the uniformizers in $R'$ tend to zero in $R^s$ with respect to the valuation topology as $K'$ ranges through a cofinal tower
 of finite extensions of $K$ in $K^t$).
\end{proof}
\fi

\sss By \eqref{sss-diag}, the lattices we have defined are related as follows:
$$ \Omega_{\mathrm{ss}}(C)\subset \Omega_{\mathrm{sat}}(C)\subset \Omega_{\mathrm{can}}(C)\otimes_R R^s.$$ We will prove in Theorem \ref{thm-logdescr} that  Edixhoven's
 jumps and Chai's elementary divisors of the Jacobian variety $\Jac(C)$ measure precisely the difference between these lattices.

\sss The saturated
 (resp.~semi-stable) lattice is invariant under base change of $C$
to a finite extension of $K$ in $K^t$ (resp.~in $K^s$).
 If $C$ is tamely ramified and $K'$ is a finite extension of
$K$ in $K^t$ such that $C\times_K K'$ has semi-stable reduction,
then
$$\Omega_{\mathrm{ss}}(C)=\Omega_{\mathrm{sat}}(C)=\Omega_{\mathrm{can}}(C\times_K K')\otimes_{R'}R^s.$$
 However, if $C$ is wildly ramified, then the intersection in the definition
of the saturated lattice {\em never} stabilizes for large $K'$\footnote{One way to see this is to combine Theorem \ref{thm-logdescr} and Corollary \ref{cor-lcm} with \cite[III.2.2.4]{HaNi-mem}: if $C$ is wildly ramified then the least common multiple of the denominators of the jumps of $\Jac(C)$ is divisible by $p$, which implies that the saturated lattice is not defined over a tame extension of $R$.}
makes the definition difficult to work with; it is not even clear from this definition that the saturated lattice is indeed a lattice.
 We will now give a
more convenient description, which also explains our choice of terminology.

\subsection{Saturated models}

\sss We will need a few results on saturated morphisms of log schemes that have been established by T.~Tsuji in an unpublished 1997 preprint; a published account of the properties
  we need can be found
 at the beginning of \cite[\S1.3]{vidal}. A morphism of saturated monoids $P\to Q$ is called saturated if, for every morphism of saturated monoids $P\to P'$, the coproduct
 $P'\oplus_P Q$ is still saturated. A morphism of $fs$ log schemes $f:X\to Y$ is called saturated if for every geometric point $x$ on $X$, the
  morphism of characteristic monoids $\overline{\mathcal{M}}_{Y,f(x)}\to \overline{\mathcal{M}}_{X,x}$ is saturated. We will only use this notion for morphisms of the form
  $f:\cC^+\to S^+$, where $\cC$ is an $R$-model of $C$. If $f$ is saturated, then for every finite separable extension $K'$ of $K$, the $fs$ base change
  $\cD^+=\cC^+\times^{fs}_{S^+}(S')^+$ coincides with the base change in the category of log schemes. In particular, the underlying scheme of $\cD^+$ is simply the fiber product
  $\cC\times_R R'$.

\sss Let $\cC$ be an $R$-model of $C$. The {\em saturation index} of $f:\cC^+\to S^+$ is defined on page 993 of \cite[\S1.3]{vidal}. It is a positive integer $m$ such that,
 for every finite separable extension $K'$ of $K$ of degree $m$, the morphism
$$\cD^+=\cC^+\times^{fs}_{S^+}(S')^+\to (S')^+$$ is saturated; we have $m=1$ if and only if $f$ is itself saturated. The saturation index $m$ is easy to compute if $\cC^+$ is log regular, which is the only case we will need: it is precisely the least common multiple of the multiplicities of the prime components in the divisor $\cC_k$.

\begin{lemma}\label{lemm-coker}
Let $P$ be a fine and saturated  monoid, and let $(\N,+)\to P$ be a
morphism of monoids. Then for every integer $d>0$, the image of
the natural morphism
$$Q=P\oplus_{\N}(1/d)\N\to Q^{\mathrm{sat}}$$ contains
$(0,1)+Q^{\mathrm{sat}}$.
\end{lemma}
\begin{proof}
 For every rational number $a$, we write $\lfloor a \rfloor$ for its integral part (the largest integer smaller than or equal to $a$) and
 $\{a\}=a-\lfloor a \rfloor$ for its fractional part. We denote by $e$ the image of $1$ under the morphism $\N\to P$.

 Using criterion (iv) in \cite[4.1]{kato-log}, it is straightforward to verify that the morphism $\N\to (1/d)\N$ is
 integral. Thus $Q$ is integral, and we can view it as the submonoid
 of the amalgamated sum $P^{\mathrm{gp}}\oplus_{\Z}(1/d)\Z$ generated by $P$ and $(0,1/d)$.

  Let $q$ be an element of $Q^{\mathrm{sat}}$. Then we can write $q$
as $(p_1,n_1/d)-(p_2,n_2/d)$ with $p_1,p_2$ in $P$ and $n_1,n_2$
non-negative integers, and we know that there exists an integer
$N>0$ such that
\begin{eqnarray*}& &N(p_1,n_1/d)-N(p_2,n_2/d)\\&=&(N(p_1-p_2)+\lfloor N(n_1-n_2)/d\rfloor
e,\{N(n_1-n_2)/d\})\end{eqnarray*}
lies in $Q$. This is only possible
 if $$N(p_1-p_2)+\lfloor N(n_1-n_2)/d\rfloor
e $$ lies in $P$, which implies that $$N(p_1-p_2+(\lfloor(n_1-n_2)/d\rfloor+1)
e) $$ lies in $P$ because $$N\lfloor n/d\rfloor +N\geq \lfloor nN/d\rfloor$$ for all integers $n$.
%
 Since $P$ is
saturated, it follows that
$$p_1-p_2+(\lfloor(n_1-n_2)/d\rfloor+1)
e$$
belongs to $P$. Thus
\begin{eqnarray*}
q+(0,1)&=&q+(e,0)\\&=&(p_1-p_2+(\lfloor (n_1-n_2)/d\rfloor +1)
e,\{(n_1-n_2)/d\})\end{eqnarray*}
lies in $Q$.
\end{proof}

\begin{lemma}\label{lemm-bounded}
We have
$$\frak{m}\Omega_{\mathrm{can}}(C)\subset \Omega_{\mathrm{log}}(C)$$ where $\frak{m}$ denotes the maximal ideal in $R$.
\end{lemma}
\begin{proof}
 Let $\cC$ be an $sncd$-model of $C$. We denote by $D$ the divisor $\cC_k- \cC_{k,\red}$ on $\cC$ and by $i$ the closed immersion $D\to \cC$.
 By Proposition \ref{prop-compar}, we have a short exact sequence of coherent $\mathcal{O}_\cC$-modules
$$0\to \clog_{\cC/R}\to \omega_{\cC/R}\to i_*i^*\omega_{\cC/R}\to 0.$$
Thus the cokernel of the inclusion of $R$-lattices
$$\Omega_{\mathrm{log}}(C)\to \Omega_{\mathrm{can}}(C)$$ is a submodule of $H^0(D,i^*\omega_{\cC/R})$. It is killed by $\frak{m}$, since every element of $\frak{m}$ vanishes on $D$.
\end{proof}

\begin{theorem}\label{thm-logtame}
Let $\cC$ be an $R$-model of $C$ such that $\cC^+$ is log regular. Let
$K'$ be a finite extension of $K$ in $K^s$ with valuation ring $R'$ such
that the log scheme
$$\cD^+=\cC^+\times^{fs}_{S^+} (S')^+$$ is saturated over $(S')^+$, and denote by $h$ the morphism
$\cD\to \cC$. Then
\begin{equation}\label{eq-sat}
\Omega_{\mathrm{sat}}(C)=H^0(\cD,h^*\clog_{\cC/R})\otimes_{R'} R^s.
\end{equation} In particular, $\Omega_{\mathrm{sat}}(C)$ is a
lattice in $H^0(C,\omega_{C/K})\otimes_K K^s$.
\end{theorem}
\begin{proof}
 Denote by $K_0$ the tame closure of $K$ in $K'$, by $R_0$ its
valuation ring, and by $S_0^+$ the spectrum of $R_0$ with its
standard log structure.
Set $C_0=C\times_K K_0$ and $\cC_0^+=\cC^+\times^{fs}_{S^+} S_0^+$ and denote by $\frak{m}_0$ the
maximal ideal of $R_0$.

 It follows from \eqref{sss-ramlog} and Lemma \ref{lemm-bounded}  that
 $$(\frak{m}_0\Omega_{\mathrm{can}}(C_0))\otimes_{R_0} R^s\subset \Omega_{\mathrm{log}}(C_0)\otimes_{R_0} R^s\subset H^0(\cD,h^*\clog_{\cC/R})\otimes_{R'}
R^s.$$
 The right hand side of this expression does not
 change if we replace $K'$ by a finite extension of $K'$ in $K^s$,
 because $\cD^+$ is saturated over $(S')^+$ so that $fs$ base change coincides with base change in the category of log
 schemes and commutes with the forgetful functor to the category of schemes.
   Since we can dominate any finite extension of $K$ in
 $K^t$ by a finite extension of $K'$ in $K^s$, and the maximal ideal $\frak{m}^s$ of $R^s$ is generated by the uniformizers
  in the finite extensions of $K$ in $K^t$,
  we find that
$$\frak{m}^s\Omega_{\mathrm{sat}}(C)\subset H^0(\cD,h^*\clog_{\cC/R})\otimes_{R'}
R^s.$$ This is only possible if $$\Omega_{\mathrm{sat}}(C)\subset H^0(\cD,h^*\clog_{\cC/R})\otimes_{R'}
R^s.$$

 Now we prove the converse inclusion.
 We claim that
\begin{equation}\label{eq-claim}
\frak{m}_0
H^0(\cD,h^*\clog_{\cC/R})\subset
H^0(\cC_0,\clog_{\cC_0/R_0})\otimes_{R_0}R'.
\end{equation}
 This implies that
$$(\frak{m}_0
H^0(\cD,h^*\clog_{\cC/R}))\otimes_{R'}R^s\subset
\Omega_{\mathrm{sat}}(C)$$ because we have
$$\Omega_{\log}(C\times_{K}K_0)\otimes_{R_0}R^s\subset \Omega_{\log}(C\times_{K}K'_0)\otimes_{R'_0}R^s\subset \Omega_{\mathrm{can}}(C\times_{K}K'_0)\otimes_{R'_0}R^s$$
 for every finite extension $K'_0$ of $K_0$ in $K^t$ with valuation ring $R'_0$.
 As in the first part of the proof, we see by replacing $K'$ by
 its finite tame extensions (and letting $K_0$ grow accordingly) that
 $$\frak{m}^s(H^0(\cD,h^*\clog_{\cD/R})\otimes_{R'}
R^s)\subset \Omega_{\mathrm{sat}}(C),$$
 and thus $$H^0(\cD,h^*\clog_{\cD/R})\otimes_{R'}
R^s\subset \Omega_{\mathrm{sat}}(C),$$ which is what we wanted to
 show.

 It remains to prove our claim \eqref{eq-claim}. The schemes $\cD$ and $\cC_0$ are
 related by a finite morphism
 $$g:\cD\to \cC_0\times_{R_0} R',$$ which is an isomorphism on
 the generic fibers. Moreover,
 $h^*\clog_{\cC/R}$ is isomorphic to the pullback of the line bundle $\clog_{\cC_0/R_0}$ to
  $\cD$. Thus it suffices to show that the cokernel of
 the morphism $$\mathcal{O}_{\cC_0\times_{R_0} R'}\to g_*\mathcal{O}_{\cD}$$ is
 killed by $\frak{m}_0$. This property is local with
 respect to the \'etale topology, so that we can assume that the
 morphism of log schemes
 $\cC^+_0\to S_0^+$ has a chart $\N\to P$, with $P$ a fine and saturated
 monoid, sending $1\in \N$ to a uniformizer $\pi_0$ in $R_0$. If we denote by $d$ the degree of $K'$ over $K_0$, then
  the morphism $(S')^+\to S_0^+$ has a chart of the form
  $$\N\to \frac{1}{d}\N\oplus \Z$$ sending $1\in \N$ to $\pi_0$, since we can write $\pi_0$ as the $d$-th power of a uniformizer in $R'$ times a  unit $u$ in $R'$.
  Thus the fiber product
 $\cC^+_0\times_{S^+_0} (S')^+$ in the category of log schemes  has a chart
 $$Q=(P\oplus \Z)\oplus_{\N}\frac{1}{d}\N\to \mathcal{O}(\cC_0\times_{S_0}S')$$
  that sends $1\in \N$ to $u^{-1}\pi_0$.
    The scheme
 $\cD$
  is given by
  $$(\cC_0\times_{R_0}R')\times_{\Z[Q]}\Z[Q^{\mathrm{sat}}],$$
 so that our claim follows from Lemma \ref{lemm-coker}.
\end{proof}

\subsection{Relation with Edixhoven's filtration and Chai's elementary divisors}
\begin{theorem}\label{thm-logdescr}\item
\begin{enumerate}
\item\label{it:relat1} The tuple of jumps in Edixhoven's filtration for $\Jac(C)$
is equal to the tuple of elementary divisors of the lattices
$$\Omega_{\mathrm{sat}}(C)\subset\Omega_{\mathrm{can}}(C)\otimes_R R^s
$$ in $H^0(C,\omega_{C/K})\otimes_K K^s$.
In particular, all the jumps are rational numbers and
$$c_{\mathrm{tame}}(\Jac(C))=c(\Omega_{\mathrm{can}}(C)\otimes_R
R^s/\Omega_{\mathrm{sat}}(C)).$$
\item\label{it:relat2} The tuple of elementary divisors for
$\Jac(C)$ is equal to the tuple of elementary divisors of the
lattices
$$
\Omega_{\mathrm{ss}}(C)\subset \Omega_{\mathrm{can}}(C)\otimes_R
R^s$$ in $H^0(C,\omega_{C/K})\otimes_K K^s$. In particular,
$$c(\Jac(C))=c(\Omega_{\mathrm{can}}(C)\otimes_R
R^s/\Omega_{\mathrm{ss}}(C)).$$
\end{enumerate}
\end{theorem}
\begin{proof}
 Point \eqref{it:relat2} is essentially the definition of the elementary divisors, modulo the identification
 $\Omega(J)=\Omega_{\mathrm{can}}(C)$ in Proposition \ref{prop-compjac} and the duality between $\Omega(J)$ and $\Lie(\mathscr{J})$.
  The same arguments show that the tuple of $K'$-elementary divisors of $\Jac(C)$ is equal to the tuple of elementary divisors of the
lattices
$$
\Omega_{\mathrm{can}}(C\times_K K')\subset \Omega_{\mathrm{can}}(C)\otimes_R
R'$$
  for every finite extension $K'$ of $K$ in $K^t$ with valuation ring $R'$.
   If we denote by
 $\frak{m}'$ the maximal ideal of $R'$, then it follows from
  \eqref{sss-diag} and  Lemma \ref{lemm-bounded} that
  $$\frak{m}'\Omega_{\mathrm{can}}(C\times_K K')\otimes_{R'}R^s\subset \Omega_{\mathrm{sat}}(C)\subset \Omega_{\mathrm{can}}(C\times_K K')\otimes_{R'}R^s.$$
  Now we apply  Proposition \ref{prop-elementary} to the chain of lattices
  $$\Omega_{\mathrm{sat}}(C)\subset \Omega_{\mathrm{can}}(C\times_K K')\otimes_{R'}R^s\subset \Omega_{\mathrm{can}}(C)\otimes_{R}R^s.$$
    This yields
    $$\frac{c_i(\Jac(C),K')}{[K':K]} \leq c_i(\Omega_{\mathrm{can}}(C)\otimes_R
R^s/\Omega_{\mathrm{sat}}(C))\leq \frac{c_i(\Jac(C),K')+1}{[K':K]}$$ for every $i$ in $\{1,\ldots,g\}$.
 Letting $K'$ range through the finite extensions of $K$ in $K^t$,
 we find that $$j_i(\Jac(C))=c_i(\Omega_{\mathrm{can}}(C)\otimes_R
R^s/\Omega_{\mathrm{sat}}(C))$$ for every $i$.
\end{proof}
 \begin{remark} Theorem \ref{thm-logdescr} suggests to define a new set of invariants, the {\em wild} elementary divisors, as the
 tuple of elementary divisors of the
lattices
$$\Omega_{\mathrm{ss}}(C)\subset \Omega_{\mathrm{sat}}(C)$$ in $H^0(C,\omega_{C/K})\otimes_K K^s$. If we define the wild base change conductor $c_{\mathrm{wild}}(Jac(C))$ as the
sum of the wild elementary divisors, then
$$c_{\mathrm{wild}}(Jac(C))=c(\Omega_{\mathrm{sat}}(C)/\Omega_{\mathrm{ss}}(C)) =c(\Jac(C))-c_{\mathrm{tame}}(\Jac(C)).$$
 The wild elementary divisors and the wild base change conductor form an interesting measure for the wild ramification of $C$; if $C$ is tamely ramified then these invariants all vanish. If $C$ is an elliptic curve, then $c_{\mathrm{wild}}(C)$ can be computed from the Swan conductor of $C$; see \cite[Ch.5,\,2.2.4]{HaNi-mem}. \end{remark}

\sss  Observe that Theorem \ref{thm-logdescr} yields a new proof of the rationality of Edixhoven's jumps for Jacobians. We will now give a more conceptual explanation for the role of the stabilization index. In order to do this, we first prove an elementary lemma.


\begin{lemma}\label{lemma-mult}
Let $\cC$ be an $sncd$-model of $C$. Let $E_0$ be an irreducible component of $\cC_k$ of multiplicity $N_0>1$, and suppose that $E_0$ is a
rational curve that meets the other components of $\cC_k$ in precisely one point. Then there exists a principal component in $\cC_k$ whose multiplicity
 $N$ is divisible by $N_0$. If $\cC$ is a minimal $sncd$-model, then we can find such a component with the additional property $N>N_0$.
 \end{lemma}
 \begin{proof}
  Let $E_1$ be the unique component of $\cC_k$ intersecting $E_0$. From the intersection formula $\cC_k\cdot E_0=0$ we obtain
  that $N_1=-N_0(E_0\cdot E_0)$; in particular, $N_0$ divides $N_1$. If $E_1$ is not principal, then it is a rational curve that meets exactly one irreducible component
  $E_2$ of $\cC_k$ different from $E_0$ (it has to meet another component because, by our assumption that $C$ has a zero divisor of degree one, the greatest common divisor of the
  multiplicities of the components of $\cC_k$ must be equal to one; see \eqref{sss-cohflat}). The intersection formula $\cC_k\cdot E_1=0$ tells us that
  $N_2=-N_0-N_1(E_1\cdot E_1)$, so that $N_0$ divides $N_2$. Repeating the argument, we eventually find a principal component $E_t$ of $\cC_k$ whose multiplicity $N_t$ is divisible by $N_0$.
   If $\cC$ is minimal, then every rational component intersecting the rest of the special fiber in at most two points has self-intersection number at most $-2$,
   because otherwise it would be contractible by Castelnuovo's criterion, contradicting the minimality of $\cC$. The above computations now easily yield $N_t>N_0$.
  \end{proof}

\begin{theorem}\label{thm-rat} If we denote by $e(C)$ the stabilization index of the curve $C$, then for every jump $j$ in Edixhoven's filtration for
 $\Jac(C)$, the product $e(C)\cdot j$ is an integer.
\end{theorem}
\begin{proof}
Let $\cC$ be an $R$-model of $C$ such that $\cC^+$ is log regular, and denote by $m$ the saturation index of
the morphism $\cC^+\to S^+$.
 By Theorem \ref{thm-logtame}, the lattice $\Omega_{\mathrm{sat}}(C)$ is defined over an extension of $K$ of degree $m$.
  Thus it follows from Theorem \ref{thm-logdescr} that the product $m\cdot j$ is integer for every jump $j$ of $\Jac(C)$. Therefore, it suffices to show that
  we can choose $\cC$ such that $m$ is equal to the stabilization index $e(C)$.

 Let $\cC'$ be the minimal $sncd$-model of $C$.  By Lipman's generalization of Artin's contractibility criterion \cite[27.1]{lipman}, any chain of rational curves in $\cC'_k$ can be contracted to a rational singularity.
     In particular, there exists a morphism  $h:\cC'\to \cC$ of normal $R$-models of $C$ that contracts precisely
  the rational components of $\cC'_k$ that meet the rest of the special fiber in exactly two points. The special fiber of $\cC$ has \'etale locally two distinct branches at any point in the image of the exceptional locus of $h$, so we can deduce from
\cite[\S3]{schroer} that $\cC^+$ is log regular.
    By Lemma \ref{lemma-mult}, the saturation index of $\cC^+\to S^+$ is equal to $e(C)$, because the non-principal components of $\cC'_k$ are either contracted by $h$ or
    do not contribute to the saturation index. This concludes the proof.
\end{proof}

\if false
\sss Unfortunately our proof of Theorem \ref{thm-rat} does not allow to conclude that $e(C)$ is the least common multiple of the denominators of the jumps of $\Jac(C)$; this was proven by the last two authors in Chapter 5, Corollary 3.1.5 of \cite{HaNi-mem} using reduction to the equicharacteristic zero case as explained in the introduction. Our proof of Theorem \ref{thm-rat} does not rely on such a reduction, but it has a ``characteristic-free flavour'', in the following sense.  Let $\cC$ be an $R$-model of $C$ such that $\cC^+$ is log regular and let $m$ be the saturation index of $\cC^+\to S^+$. Let $K'$ be a separable degree $m$  extension of $K$ and set $\cD^+=\cC^+\times^{fs}_{S^+}(S')^+$. If $m$ is prime to $p$ (in particular, if $k$ has characteristic zero) then $\cD$ is precisely the normalization of $\cC\times_R R'$. If $p>1$ and $p$ divides $m$ then $\cD$ is only a partial normalization of $\cC\times_R R'$.
 Its construction relies on the combinatorics of the log structure on $\cC^+$ and is less sensitive to the arithmetic properties of $\cC$  than the full normalization of $\cC\times_R R'$, which depends strongly on the choice of the wild extension $K'$.
\fi

\section{A formula for the jumps}
\subsection{The basic formula}
\sss \label{sss-formula} We will now use our logarithmic interpretation of the jumps of $\Jac(C)$ in Theorem \ref{thm-logdescr} to deduce an explicit formula for the jumps
 in terms of the combinatorial
reduction data of $C$.
 Our starting point is the chain of lattices
$$\Omega_{\log}(C)\otimes_R R^s\subset \Omega_{\mathrm{sat}}(C)\subset  \Omega_{\mathrm{can}}(C)\otimes_R R^s$$ in $H^0(C,\omega_{C/K})\otimes_K K^s$.
 Recall that, by Theorem \ref{thm-logdescr}, the jumps of $\Jac(C)$ are precisely the elementary divisors
 of the pair of lattices $$\Omega_{\mathrm{sat}}(C)\subset  \Omega_{\mathrm{can}}(C)\otimes_R R^s.$$
  We will compute these elementary divisors by computing those of the other inclusions in the chain:
  \begin{eqnarray}
  \Omega_{\log}(C) &\subset& \Omega_{\mathrm{can}}(C), \label{eq-incl1}
  \\ \Omega_{\log}(C)\otimes_R R^s &\subset& \Omega_{\mathrm{sat}}(C). \label{eq-incl2}
  \end{eqnarray}
    The tuples of elementary divisors do not behave additively in chains, in general, but they do in special cases, as is explained in  the following
   easy lemma.

   \begin{lemma}\label{lemm-eldiv}
   Let $V$ be a vector space over $K$ of finite dimension $g$ and let $\Omega_1\subset \Omega_2\subset \Omega_3$ be a chain of lattices
   in $V\otimes_K K^s$. Suppose that the tuple of elementary divisors of $\Omega_1\subset \Omega_3$ is of the form
   $$v=(0,\ldots,0,a,\ldots,a)$$ for some positive rational number $a$, and denote by $n\leq g$ the number of entries equal to $0$. Then the tuple of elementary divisors
    $$w=(w_1,\ldots,w_g)$$ of $\Omega_1\subset \Omega_2$ satisfies $0\leq w_i\leq v_i$ for all $i$ in $\{1,\ldots,g\}$, and the tuple of elementary divisors
    of $\Omega_2\subset \Omega_3$ is given by
    $$(0,\ldots,0,a-w_{g},\ldots,a-w_{n+1}).$$
   \end{lemma}
   \begin{proof}
    We may assume that all the lattices are defined over $R$. Then $a$ is a positive integer and the statement boils down to the following simple fact:
    if $N$ is a submodule of $M=(R/\frak{m}^a)^{g-n}$ for some positive integer $a$, then $N$ is isomorphic to
    $\oplus_{i=1}^{g-n}R/\frak{m}^{a_i}$ for some non-negative integers $a_i\leq a$, and the quotient $M/N$ is isomorphic to $\oplus_{i=1}^{g-n}R/\frak{m}^{a-a_i}$.
   \end{proof}

The following proposition will be quite useful in our computations.

\begin{prop}\label{prop-torsfree}
If $\cC$ is an $sncd$-model of $C$, then the $R$-modules
$H^1(\cC,\omega_{\cC/R})$ and $H^1(\cC,\clog_{\cC/R})$ have no torsion.
\end{prop}
\begin{proof} We write $D$ for the divisor $\cC_k-\cC_{k,\red}$ on $\cC$.
We know by Proposition \ref{prop-compar} that $$\clog_{\cC/R}=\omega_{\cC/R}(-D).$$
 The cohomology module $H^1(\cC,\mathcal{O}_{\cC})$ is free by cohomological flatness of the structural morphism $\cC\to S$ (see \eqref{sss-cohflat}).
 The short exact sequence of $\mathcal{O}_{\cC}$-modules
$$0\to \mathcal{O}_{\cC}(D-\cC_k)\to \mathcal{O}_{\cC}\to \mathcal{O}_{\cC_{k,\red}}\to 0$$ gives rise to an injection
  $$H^1(\cC,\mathcal{O}_{\cC}(D))\to H^1(\cC,\mathcal{O}_{\cC})$$ by surjectivity of the map
  $$R=H^0(\cC,\mathcal{O}_{\cC})\to H^0(\cC,\mathcal{O}_{\cC_{k,\red}})=k.$$ Thus we see that
  $H^1(\cC,\mathcal{O}_{\cC}(D))$ is free, as well.
 Now Grothendieck-Serre duality provides us with isomorphisms
$$H^1(\cC,\omega_{\cC/R})\to H^0(\cC,\mathcal{O}_{\cC})^{\vee} $$
 and
 $$H^1(\cC,\clog_{\cC/R})\to H^0(\cC,\mathcal{O}_{\cC}(D))^{\vee} .$$ In particular, $H^1(\cC,\omega_{\cC/R})$ and $H^1(\cC,\clog_{\cC/R})$ have no torsion.
\end{proof}


\subsection{Computation of $\Omega_{\mathrm{can}}(C)/\Omega_{\mathrm{log}}(C)$}
\sss 
 We know that $\frak{m}\Omega_{\mathrm{can}}(C)\subset \Omega_{\log}(C)$ by Lemma \ref{lemm-bounded}.
  Thus the quotient $\Omega_{\mathrm{can}}(C)/\Omega_{\log}(C)$ is isomorphic to $(R/\frak{m})^u$ for some  $0\leq u\leq g$ and the tuple of elementary divisors of
  $\Omega_{\log}(C)\subset \Omega_{\mathrm{can}}(C)$ is given by $(0,\ldots,0,1,\ldots,1)$ where the number of zeroes is equal to $g-u$.
  It only remains to determine the value of $u$, that is, the dimension of the $k$-vector space $\Omega_{\mathrm{can}}(C)/\Omega_{\log}(C)$.
  To compute this dimension, we rewrite it as
$$\mathrm{dim}_k (\Omega_{\mathrm{can}}(C)/\Omega_{\log}(C))=g- \mathrm{dim}_k (\Omega_{\log}(C)/\frak{m}\Omega_{\mathrm{can}}(C)).$$
 Let $\cC$ be an $sncd$-model of $C$. By Proposition \ref{prop-compar}, we have a short exact sequence of $\mathcal{O}_{\cC}$-modules
\begin{equation}\label{eq-exactseq1}
0\to \frak{m}\omega_{\cC/R}=\clog_{\cC/R}(-\cC_{k,\red})\to \clog_{\cC/R}\to \iota_*\iota^* \clog_{\cC/R}\to 0\end{equation}
 where $\iota$ denotes the closed immersion $\cC_{k,\red}\to \cC$. By
 Proposition \ref{prop-torsfree}, the $R$-module $H^1(\cC,\omega_{\cC/R})$ is torsion free. Therefore, the sequence
 $$0\to H^0(\cC,\frak{m}\omega_{\cC/R})\to H^0(\cC,\clog_{\cC/R})\to H^0(\cC_{k,\red},\iota^*\clog_{\cC/R})\to 0$$ is still exact; the start of this sequence is precisely the inclusion of lattices $\frak{m}\Omega_{\mathrm{can}}(C)\to \Omega_{\log}(C)$. Hence,
$$\mathrm{dim}_k (\Omega_{\mathrm{can}}(C)/\Omega_{\log}(C))=g-\mathrm{dim}_k H^0(\cC_{k,\red},\iota^*\clog_{\cC/R}).$$

\sss \label{sss-resol} Now we compute $\mathrm{dim}_k H^0(\cC_{k,\red},\iota^*\clog_{\cC/R})$. Denote by
$$a_1:\widetilde{\cC}_{k,\red}\to  \cC_{k,\red}$$ the normalization morphism of $\cC_{k,\red}$. If we write $\cC_{k,\red}=\sum_{i=1}^r E_i$
 then $\widetilde{\cC}_{k,\red}$ is simply the disjoint union $\sqcup_i E_i$. We denote by $\Sigma$ the singular locus of
  $\cC_{k,\red}$ and by $a_2$ the closed immersion $\Sigma\to \cC_{k,\red}$. For notational convenience, we write $\mathcal{L}$  for the line bundle  $\clog_{\cC/R}$ on $\cC$, and $\mathcal{L}_i$ for the pullback of $\mathcal{L}$ to $E_i$.
   Then we can construct the usual short exact sequence
  $$0\to \iota^*\mathcal{L} \to (a_1)_*(a_1)^*\iota^*\mathcal{L}\to (a_2)_*(a_2)^*\iota^*\mathcal{L}\to 0$$
   which tells us that $$\chi(\iota^*\mathcal{L})=\sum_{i=1}^r \chi(\mathcal{L}_i) -|\Sigma|.$$
     By Proposition \ref{prop-compar} and the adjunction formula, $\mathcal{L}_i$ is isomorphic to
   the sheaf $$\omega_{E_i/k}(E_i\cap \Sigma)$$ of differential forms on $E_i$ with logarithmic poles along $E_i\cap \Sigma$. Thus
 $$\chi(\mathcal{L}_i)=g(E_i)-1+|E_i\cap \Sigma|$$ where $g(E_i)$ denotes the genus of $E_i$. Since every point of $\Sigma$ lies on precisely two irreducible components of
 $\cC_{k,\red}$, we find
 $$\chi(\iota^*\mathcal{L})=\sum_{i=1}^r g(E_i)-r +|\Sigma|$$
 and hence
 $$\mathrm{dim}_k H^0(\cC_{k,\red},\iota^*\mathcal{L})= \sum_{i=1}^r g(E_i)-r +|\Sigma| + \mathrm{dim}_k H^1(\cC_{k,\red},\iota^*\mathcal{L}).$$

 \sss   The short exact sequence \eqref{eq-exactseq1} also yields an exact sequence
  $$0\to H^1(\cC,\clog_{\cC/R}(-\cC_{k,\red}))\to H^1(\cC,\clog_{\cC/R}) \to H^1(\cC_{k,\red},\iota^*\mathcal{L})\to 0.$$
   By Proposition \ref{prop-compar} and Grothendieck-Serre duality, the second arrow can be identified with the map
   $$H^0(\cC,\mathcal{O}_{\cC}(\cC_k))^{\vee}\to H^0(\cC,\mathcal{O}_{\cC}(\cC_k- \cC_{k,\red}))^{\vee}$$
   whose cokernel is clearly isomorphic to $k$. Thus $$\mathrm{dim}_k H^1(\cC_{k,\red},\iota^*\mathcal{L})=1.$$
   Putting all these pieces together, we arrive at the following result.

   \begin{prop}\label{prop-term1}
   Let $\cC$ be an $sncd$-model of $C$. Denote by $\Gamma$ the dual graph of $\cC_k$ and by $\beta(\Gamma)$ its first Betti number. We write $u(C)$ for the unipotent rank of
   $\Jac(C)$, that is, the dimension of the unipotent radical of the identity component of the special fiber of the N\'eron model of $\Jac(C)$. Then
  \begin{eqnarray*}
  \mathrm{dim}_k (\Omega_{\mathrm{can}}(C)/\Omega_{\log}(C))&=&g-\sum_{i=1}^r g(E_i)-\beta(\Gamma)
  \\ &=&u(C).
  \end{eqnarray*}
  \end{prop}
  \begin{proof}
  The above computations yield the formula
     $$\mathrm{dim}_k (\Omega_{\mathrm{can}}(C)/\Omega_{\log}(C))=g-\sum_{i=1}^r g(E_i)+r-|\Sigma|-1.$$
 Since $r$ is equal to the number of vertices of $\Gamma$ and $|\Sigma|$ equals the number of edges, we have
  $r-|\Sigma|-1=-\beta(\Gamma)$. The equality  $$g-\sum_{i=1}^r g(E_i)-\beta(\Gamma)=u(C)$$ is well-known; see for instance \cite[p.148]{lorenzini-comp}.
  \end{proof}

\begin{remark} If $\cC_k$ is reduced, then $\Omega_{\mathrm{can}}(C)=\Omega_{\log}(C)$ and our formula in Proposition \ref{prop-term1} boils down
   to the classical expression
   $$g= \sum_{i=1}^r g(E_i)+\beta(\Gamma).$$
\end{remark}

\subsection{Computation of $\Omega_{\mathrm{sat}}(C)/(\Omega_{\log}(C)\otimes_R R^s )$}
\sss We still denote by $\cC$ an $sncd$-model of $C$ and by $\mathcal{L}$ the line bundle $\clog_{\cC/R}$ on $\cC$.
 Denote by $m$ the saturation index  of $\cC^+\to S^+$, that is, the least
 common multiple of the multiplicities of the prime components in $\cC_k$.
 Let $K'$ be
 a degree $m$ extension of $K$ in $K^s$ and denote by $R'$ its valuation ring. Set $\cC'=\cC\times_R R'$ and denote by $\cD$ the underlying scheme of $\cC^+\times^{fs}_{S^+}(S')^+$.
 We write $\mathcal{L}'$ for the pullback of $\mathcal{L}$ to $\cC'$ and $h$ for the natural morphism $\cD\to \cC'$. By flat base change, we have $$H^0(\cC',\mathcal{L}')=\Omega_{\log}(C)\otimes_R R'.$$
  By Theorem \ref{thm-logtame}, we know that
 $$\Omega_{\mathrm{sat}}(C)=H^0(\cD,h^*\mathcal{L}')\otimes_{R'}R^s.$$ Therefore, we need to compute the cokernel of
 $$H^0(\cC',\mathcal{L}')\to H^0(\cD,h^*\mathcal{L}').$$

\sss Denote by $\mathcal{F}$ the cokernel of the morphism of $\mathcal{O}_{\cC'}$-modules
 $$\mathcal{O}_{\cC'}\to h_*\mathcal{O}_\cD.$$ Note that $\mathcal{F}$ is trivial on the generic fiber $\cC'$ because the morphism
 $\cD_{K'}\to \cC'_{K'}$ is an isomorphism.
   We have a short exact sequence of $\mathcal{O}_{\cC'}$-modules
$$0\to \mathcal{L}'\to h_*h^*\mathcal{L}'\to \mathcal{L}'\otimes_{\mathcal{O}_{\cC'}} \mathcal{F} \to 0$$ which gives rise to a short exact sequence of $R'$-modules
$$0\to H^0(\cC,\mathcal{L}')\to H^0(\cD,h^*\mathcal{L}')\to H^0(\cC',\mathcal{L}'\otimes_{\mathcal{O}_{\cC'}} \mathcal{F})\to 0 $$
because $H^1(\cC',\mathcal{L}')=H^1(\cC,\mathcal{L})\otimes_R R'$ has no torsion (see Proposition \ref{prop-torsfree}).
 As we explained in the proof of Theorem \ref{thm-logtame}, it follows from Lemma \ref{lemm-coker} that $\mathcal{F}$ is
 killed by $\frak{m}$, so that we can also view it as a coherent sheaf on the $k$-scheme $$X=\cC'\times_{R'}(R'/\frak{m}R')\cong \cC_k\times_k (R'/\frak{m}R').$$ Thus, denoting by $\mathcal{L}_k$ the pullback of $\mathcal{L}$ to $\cC_k$ and by $\kappa$ the projection morphism $X\to \cC_k$, we find that the cokernel of
 $$H^0(\cC',\mathcal{L}')\to H^0(\cD,h^*\mathcal{L}')$$
 is isomorphic to the $R'/\frak{m}R'$-module
 $$M=H^0(\cC_k,\mathcal{L}_k\otimes_{\mathcal{O}_{\cC_k}} \kappa_*\mathcal{F}).$$

\sss \label{sss-filtr} Denote by $\frak{m}'$ the maximal ideal of $R'$.
For every $i$ in $\{0,\ldots,m\}$, we set
$$\mathcal{F}_i=\kappa_*((\frak{m}')^i\mathcal{F}/(\frak{m}')^{i+1}\mathcal{F})$$
 and $V_i=(\frak{m}')^iM/(\frak{m}')^{i+1}M$.
Then $\mathcal{F}_i$ is a coherent $\mathcal{O}_{\cC_k}$-module and $V_i$ is a $k$-vector space, and they both vanish if $i=m$. The dimensions of the vector spaces $V_i$ completely determine the $R'$-module structure of $M$: for every $i$ in $\{0,\ldots,m-1\}$, the multiplicity of $R'/(\frak{m}')^{i+1}$ as a direct summand in $M$
   is equal to $\mathrm{dim}_kV_{i}-\mathrm{dim}_kV_{i+1}$.

 From the spectral sequence for the hypercohomology of a filtered complex \cite[1.4.5]{deligne-hodgeII}, we deduce that
\begin{equation}\label{eq-connect}
V_i\cong \ker(d_i:H^0(\cC_k,\mathcal{L}_k\otimes\mathcal{F}_i)\to H^1(\cC_k,\mathcal{L}_k\otimes\mathcal{F}_{i+1}) )\end{equation}
  for every  $i$ in $\{0,\ldots,m-1\}$, where $d_i$ is the connecting homomorphism in the long exact cohomology sequence associated to the short exact sequence of
   $\mathcal{O}_{\cC_k}$-modules $$ 0\to \mathcal{L}_k\otimes\mathcal{F}_{i+1}\to
   \mathcal{L}_k\otimes \kappa_*((\frak{m}')^{i}\mathcal{F}/(\frak{m}')^{i+2}\mathcal{F})
    \to \mathcal{L}_k\otimes \mathcal{F}_i\to 0.  $$

\sss \label{sss-divnot} We will now give an explicit description of the $\mathcal{O}_{\cC_k}$-modules $\mathcal{F}_i$. It will be convenient to use the
  following notation. For every effective $\R$-divisor $D$ on $\cC$, we denote by $\lfloor D \rfloor$ the integral part of $D$, obtained by rounding down the coefficients, and we
   write $\langle D\rangle $ for the divisor obtained from $D$ by putting all the coefficients in $\R\setminus \Z$ equal to zero.
     Moreover, we write $\langle D\rangle_{\red}$ for the reduction of $\langle D\rangle $ and $\mathcal{J}(D)$ for the pullback of $\mathcal{O}_{\cC}(\lfloor D\rfloor)$ to $\langle D\rangle_{\red}$.
  The letter $\mathcal{J}$ stands for ``jump'': if we let the coefficients of $D$ grow continuously in $\R$, then $\mathcal{J}(D)$ detects how $\lfloor D\rfloor$ changes.

\begin{prop}\label{prop-split}
For every $i$ in $\{0,\ldots,m-1\}$, the $\mathcal{O}_{\cC_k}$-module $\mathcal{F}_i$ is isomorphic to
$$\bigoplus_{j=1}^{m-i-1} \mathcal{J}((j/m)\cC_{k}).$$
\end{prop}
\begin{proof}
The proof is based on a local study of the morphism $h:\cD\to \cC'$.

{\em Case 1}. First, let $x$ be a regular point of $\cC_{k,\red}\cong \cC'_{k,\red}$.
  Since $\cC$ is an $sncd$-model, the morphism $\cC^+\to S^+$ has Zariski-locally at $x$ a chart of the form $$\N\to \frac{1}{a}\N\oplus \Z:1\mapsto (1,1)$$
   where $a$ is the multiplicity of $\cC_k$ at $x$, the generator of  $\N$ is sent to a uniformizer $\pi$ in $R$, and the generator of $(1/a)\N$ is sent to a local defining function $f$ for $\cC_{k,\red}$ in $\cC$ at $x$.
  As we've already explained in the proof of Theorem \ref{thm-logtame}, the morphism $(S')^+\to S^+$ has a chart of the form $$\N\to \frac{1}{m}\N\oplus \Z:1\mapsto (1,1)$$ sending the generator of $\N$ to $\pi$ and the generator of $(1/m)\N$ to a uniformizer $\pi'$ in $R'$. Thus, locally at $x$, $(\cC')^+$ has a chart of the form
  $$Q=(\Z\oplus\frac{1}{a}\N)\oplus_{\N}\frac{1}{m}\N\to \mathcal{O}_{\cC',x}$$ sending  the generator of $(1/a)\N$ to $f$ and the generator of $(1/m)\N$ to $\pi'$. If we denote by $Q\to Q^{\sat}$ the natural morphism from $Q$ to its saturation, then over some open neighbourhood of $x$ in $\cC'$, $\cD$ is given by $$\cC'\times_{\Z[Q]}\Z[Q^{\sat}].$$ Thus we must understand the exact shape of the morphism $Q\to Q^{\sat}$.

   The groupification $Q^{\gp}$ of $Q$ is the coproduct
   $$(\Z \oplus \frac{1}{a}\Z)\oplus_{\Z}\frac{1}{m}\Z\cong (\Z \oplus \frac{1}{a}\Z \oplus \frac{1}{m}\Z)/\langle (1,1,-1)\rangle.$$
  An element $(u,v/a,w/m)$ of $Q^{\gp}$ belongs to $Q^{\sat}$ if and only if $v/a+w/m$ is non-negative, that is, $mv+aw\geq 0$. Thus the $\Z[Q]$-module
  $\Z[Q^{\sat}]/\Z[Q]$ is generated by the elements $(0,-\lfloor ja/m \rfloor ,j)$ with $j\in \{1,\ldots,m-1\}$. This means that the stalk $\mathcal{F}_x$
   is generated as an $\mathcal{O}_{\cC,x}$-module by the elements $(\pi')^j/f^{\lfloor ja/m \rfloor}$ with $j\in \{1,\ldots,m-1\}$. An element of the form $(\pi')^s/f^t$ is divisible by $(\pi')^i$ in $\mathcal{F}_x$
    if and only if $a(s-i)-mt\geq 0$. In particular, if $(\pi')^s/f^t$ is divisible by $(\pi')^i$ then $(\pi')^s/f^{t-1}$ is divisible by $(\pi')^{i+1}$.
       From these computations, we deduce that locally at $x$, the morphism
   $$\bigoplus_{j=1}^{m-1-i} \mathcal{O}_{\cC}(\lfloor (j/m) \cC_k\rfloor )\to \mathcal{F}_i:(c_{1},\ldots,c_{m-1-i})\mapsto \sum_{j=1}^{m-1-i}c_j (\pi')^{j+i}$$
   factors through an isomorphism $$\bigoplus_{j=1}^{m-i-1} \mathcal{J}((j/m)\cC_{k})\to \mathcal{F}_i,$$ for every $i\in \{0,\ldots,m-1\}$.

{\em Case 2}.    Now we treat the case where $x$ is a singular point of $\cC_{k,\red}$.
 The morphism $\cC^+\to S^+$ has Zariski-locally at $x$ a chart of the form $$\N\to \Z\oplus \frac{1}{a}\N\oplus \frac{1}{b}\N:1\mapsto (1,1,1)$$
   where $a$ and $b$ are the multiplicities of $\cC_k$ along the irreducible components $E_1$ and $E_2$ passing through at $x$, the generator of  $\N$ is sent to a uniformizer $\pi$ in $R$, and the generators $(1/a)\N$ and $(1/b)\N$ are sent to  local defining functions $f_1$ and $f_2$ for $E_1$ and $E_2$ in $\cC$ at $x$. In a similar way as in Case 1, we can
   produce a local chart $$Q=(\Z\oplus\frac{1}{a}\N\oplus \frac{1}{b}\N)\oplus_{\N}\frac{1}{m}\N\to \mathcal{O}_{\cC',x}$$
   sending  the generators of $(1/a)\N$ and $(1/b)\N$ to $f_1$ and $f_2$, respectively, and the generator of $(1/m)\N$ to $\pi'$.
    Now $Q^{\gp}$ is given by $$(\Z\oplus\frac{1}{a}\Z\oplus \frac{1}{b}\Z)\oplus_{\Z}\frac{1}{m}\Z\cong (\Z\oplus \frac{1}{a}\Z\oplus \frac{1}{b}\Z \oplus \frac{1}{m}\Z)/\langle (1,1,1,-1)\rangle $$
     and an element $(t,u/a,v/b,w/m)$ of $Q^{\gp}$ belongs to $Q^{\sat}$ if and only if  $mu+aw\geq 0$ and $mv+bw\geq 0$. It follows that the stalk $\mathcal{F}_x$
   is generated over $\mathcal{O}_{\cC,x}$ by the elements $(\pi')^j/(f_1^{\lfloor ja/m \rfloor}f_2^{\lfloor jb/m \rfloor})$ with $j\in \{1,\ldots,m-1\}$ and $j\geq \min\{m/a,m/b\}$. We deduce as in Case 1 that locally at $x$, the morphism
   $$\bigoplus_{j=1}^{m-i-1} \mathcal{O}_{\cC}(\lfloor (j/m)\cC_{k} \rfloor )\to \mathcal{F}_i:(c_1,\ldots,c_{m-1-i})\mapsto \sum_{j=1}^{m-1-i}c_j (\pi')^{j+i}$$
   factors through an isomorphism $$\bigoplus_{j=1}^{m-i-1} \mathcal{J}((j/m)\cC_{k})\to \mathcal{F}_i,$$ for every $i\in \{0,\ldots,m-1\}$.

 As $x$ varies, all of these local isomorphisms glue to an isomorphism as required in the statement.
\if false
  Now assume that $m$ is strictly larger than $m_0$; since we assume that $m$ is a multiple of $m_0$, this means that $m\geq 2m_0$. Using the above computations, one sees that for every $i\in \{0,\ldots,m-1\}$
 the $\mathcal{O}_{\cC}$-module   $$\mathcal{G}_i=\kappa_*((\frak{m}')^{i}\mathcal{F}/(\frak{m}')^{i+2}\mathcal{F})$$ is still killed by the defining ideal sheaf of $\cC_{k,\red}$ in
 $\cC$. Indeed, in Case 1,  $(\pi')^j/f^{\lfloor ja/m \rfloor-1}$ is divisible by $(\pi')^2$ in $\mathcal{F}_x$ because $m\geq 2a$ so that
 $$\lfloor ja/m \rfloor -1 \leq \lfloor (j-2)a/m\rfloor . $$ Likewise, in Case 2, $(\pi')^j/(f_1^{\lfloor ja/m \rfloor-1}f_2^{\lfloor jb/m \rfloor-1})$ is divisible by $(\pi')^2$ because
  $m\geq 2a$ and $m\geq 2b$. Thus we can view $\mathcal{G}_i$ as an $\mathcal{O}_{\cC_{k,\red}}$-module, and we see that the map
  $$\bigoplus_{j=1}^{m-i-1} \mathcal{J}((j/m)\cC_{k})\to \mathcal{G}_i:(c_1,\ldots,c_{m-1-i})\mapsto \sum_{j=1}^{m-1-i}c_j (\pi')^{j+i} $$
  is a section of the projection morphism $\mathcal{G}_i \to \mathcal{F}_i$.
  Thus the sequence \eqref{eq-sequence} is split.
 \fi
\end{proof}

\sss We can now use Proposition \ref{prop-split} to finish our computations. We will need the following vanishing result.
  Recall that we write $\cC_k=\sum_{i=1}^r N_i E_i$ and that
 we denote by $\mathcal{L}$ the line bundle $\clog_{\cC/R}$ on $\cC$, by $\iota$ the closed immersion $\cC_{k,\red}\to \cC_k$, and by $\Sigma$ the
 set of singular points of $\cC_{k,\red}$.

\begin{lemma}\label{lemma-vanish}
Let $E$ be a connected, smooth and proper curve over $k$ and let $D$ and $D'$ be divisors on $E$ of degrees $d$ and $d'$, respectively. Assume that $D'$ is reduced.
Then $H^1(E,\omega_{E/k}\otimes \mathcal{O}_{E}(D))$ vanishes if $d>0$, and the restriction map  $$H^0(E,\omega_{E/k}\otimes \mathcal{O}_{E}(D))\to \bigoplus_{x\in D'}\iota_x^*(\omega_{E/k}\otimes \mathcal{O}_{E}(D))$$ is surjective if $d>d'$. Here we wrote $\iota_x$ for the closed immersion of $x$ in $E$.
\end{lemma}
\begin{proof}
These are standard applications of Serre duality: if $d>0$ then
$$H^1(E,\omega_{E/k}\otimes \mathcal{O}_{E}(D))\cong H^0(E,\mathcal{O}_{E}(-D))^{\vee}=0.$$
 Moreover, the cokernel of $$H^0(E,\omega_{E/k}\otimes \mathcal{O}_{E}(D))\to \bigoplus_{x\in D'}\iota^*_x(\omega_{E/k}\otimes \mathcal{O}_{E}(D))$$ is isomorphic to
 $$H^1(E,\omega_{E/k}\otimes \mathcal{O}_{E}(D-D')),$$
 and thus vanishes if $d>d'$. 
\end{proof}
\begin{prop}\label{prop-vanish}
For every $j$ in $\{1,\ldots,m-1\}$ we have
$$H^1(\cC_{k,\red},\iota^*\mathcal{L}\otimes\mathcal{J}((j/m)\cC_{k}))=0.$$
\end{prop}
\begin{proof}
We denote by $I_j$ the set of indices $i$ in $\{1,\ldots,r\}$ such that $j$ is a multiple of $m/N_i$.
 The components $E_i$ with $i\in I_j$ are precisely the prime components of $\langle (j/m)\cC_k\rangle$.
 For every $i\in I_j$ we denote by $\mathcal{J}_{j,i}$ the restriction of
 $\mathcal{J}((j/m)\cC_k)$ to $E_i$. This is a line bundle on $E_i$ whose degree  is equal to the intersection number $(E_i\cdot \lfloor (j/m)\cC_k\rfloor)$. Writing
 $\{(j/m)\cC_k\}$ for the fractional part $(j/m)\cC_k-\lfloor (j/m)\cC_k\rfloor$ of the $\Q$-divisor $(j/m)\cC_k$, we compute:
 $$0=(E_i \cdot (j/m)\cC_k)= (E_i\cdot \lfloor (j/m)\cC_k\rfloor) + (E_i\cdot \{(j/m)\cC_k\}).$$
 The prime components of the divisor $\{(j/m)\cC_k\}$ are precisely the components of $\cC_k$ that are not contained in $\langle (j/m)\cC_k\rangle$. We write $\sigma_{j,i}$ for the number of intersection points of $E_i$ with the support of $\{(j/m)\cC_k\}$. Note that these are singular points of $\cC_{k,\red}$, so that
    $\sigma_{j,i}\leq |\Sigma\cap E_i|$.
      Since
   the coefficients of $\{(j/m)\cC_k\}$ are strictly contained between $0$ and $1$, we find that
   $$(E_i\cdot \lfloor (j/m)\cC_k\rfloor)\geq \min\{-\sigma_{j,i} +1,0\}.$$
  We have already observed in \eqref{sss-resol} that the restriction $\mathcal{L}_i$ of $\mathcal{L}$ to $E_i$ is isomorphic to the sheaf of differential forms with logarithmic poles at the points of
  $(\Sigma\cap E_i).$  Thus Lemma \ref{lemma-vanish} implies that
  $$H^1(E_i,\mathcal{L}_i\otimes\mathcal{J}_{j,i})=0$$ for every $i\in I_j$.
     Now we can compute  $$H^1(\cC_{k,\red},\iota^*\mathcal{L}\otimes\mathcal{J}((j/m)\cC_{k}))$$ using the standard resolution of
  the line bundle $\iota^*\mathcal{L}\otimes\mathcal{J}((j/m)\cC_{k})$ on the reduced strict normal crossings divisor $\langle (j/m)\cC_k\rangle_{\red}$ as in \eqref{sss-resol}.
 The associated long exact cohomology sequence contains the exact subsequence
 $$\bigoplus_{i\in I_j}H^0(E_i,\mathcal{L}_i\otimes  \mathcal{J}_{j,i})\to \bigoplus_{x\in \Sigma_j}\iota_x^*(\mathcal{L}\otimes  \mathcal{J}((j/m)\cC_{k}))\to H^1(\cC_{k,\red},\iota^*\mathcal{L}\otimes\mathcal{J}((j/m)\cC_{k}))\to 0$$
  where $\Sigma_j$ denotes the set of singular points of $\langle (j/m)\cC_k\rangle_{\red}$ and we write $\iota_x$ for the closed immersion of $x$ in $\cC_{k,\red}$.
   We will prove that the first arrow in this sequence is surjective, which means that $$H^1(\cC_{k,\red},\iota^*\mathcal{L}\otimes\mathcal{J}((j/m)\cC_{k}))=0.$$

  We say that a point $x_0$ in $\Sigma_j$ is {\em good} if the characteristic function of $\{x_0\}$ lies in the image of
  $$\rho:\bigoplus_{i\in I_j}H^0(E_i,\mathcal{L}_i\otimes  \mathcal{J}_{j,i})\to \bigoplus_{x\in \Sigma_j}\iota_x^*(\mathcal{L}\otimes  \mathcal{J}((j/m)\cC_{k})).$$
 For every $i\in I_j$, the restriction map
 $$ H^0(E_i,\mathcal{L}_i\otimes  \mathcal{J}_{j,i})\to \bigoplus_{x\in E_i\cap \Sigma_j}\iota_x^*(\mathcal{L}_i\otimes  \mathcal{J}_{j,i})$$
 is surjective if $\sigma_{j,i}\neq 0$,
  because the number of points in $E_i\cap \Sigma_j$ is precisely equal to $|\Sigma\cap E_i|-\sigma_{j,i}$ so that we can apply Lemma \ref{lemma-vanish}.
  Thus if if $\sigma_{j,i}\neq 0$ then every point of $\Sigma_j\cap E_i$ is good. On the other hand, if $i$ is any element of $I_j$ such that $\sigma_{j,i}=0$ and $E_i$ contains a good point $x_0$ of
  $\Sigma_j$, then any point $x_1\neq x_0$ of $\Sigma_j\cap E_i$ is good, because we can always find an element of
  $$H^0(E_i,\mathcal{L}_i\otimes  \mathcal{J}_{j,i})$$ that is non-zero at $x_0$ and $x_1$ and vanishes at all other points of $\Sigma_j\cap E_i$ (this is again a consequence of
   Lemma \ref{lemma-vanish}, applied to $D'=\{x_0\}$ and $D'=\{x_1\}$).
   Note that $I_j$ cannot be the whole set $\{1,\ldots,r\}$ because this contradicts our overall assumption that the multiplicities $N_i$ are coprime. Thus every connected component of $\sum_{i\in I_j}E_i$ contains at least one prime component $E_i$ with $\sigma_{j,i}\neq 0$, and we may conclude that every point in $\Sigma_j$ is good, so that $\rho$ is surjective.
\end{proof}

\subsection{An explicit formula for the jumps}

\begin{theorem}\label{theorem-formula}
 Let $\cC$ be an $sncd$-model of $C$. Write $\cC_k=\sum_{i=1}^r N_i E_i$ and denote by $m$ the least common multiple of the multiplicities $N_i$ of the irreducible components $E_i$ in $\cC_k$. For each $i$, we denote by $g(E_i)$ the genus of $E_i$. For every $j$ in $\{0,\ldots,m\}$, we denote by $I_j$ the set of indices $i$ in $\{1,\ldots,r\}$ such that
 $j$ is a multiple of $m/N_i$, and we write $\sigma_j$ for the number of singular points of $\cC_{k,\red}$ that lie on at least one component $E_i$ with $i\in I_j$.

 The jumps of $\Jac(C)$ are contained in the set $\{0,\ldots,(m-1)/m\}$.  For every element $j$ in $\{0,\ldots,m-1\}$, the multiplicity of $j/m$ as a jump of $\Jac(C)$ is equal to
  $$(\sum_{i\in I_j}E_i\cdot \lfloor  (j/m)\cC_{k} \rfloor )+ \sum_{i\in I_j}g(E_i)-|I_j|+\sigma_j+\delta_{j,0}$$
 where  $\delta$ is the Kronecker symbol.
\end{theorem}
\begin{proof}
 We already know by Theorem \ref{thm-rat} that the jumps are multiples of $1/m$, and it follows from Edixhoven's original definition in \cite{edix} that they are always strictly smaller than one (for arbitrary abelian $K$-varieties). Thus each jump is of the form $j/m$ with $j\in \{0,\ldots,m-1\}$.
  First, assume that $j\neq 0$. By \eqref{sss-formula} and Lemma \ref{lemm-eldiv}, the multiplicity of $j/m$ as a jump of $\Jac(C)$ is equal to the
  number  of occurrences of the value $(m-j)/m$ in the tuple of elementary divisors of the inclusion of lattices
 $$\Omega_{\log}(C)\otimes_R R^s \subset \Omega_{\mathrm{sat}}(C).$$
 We will denote this number by $\gamma_{j}$.
      As we have explained in \eqref{sss-filtr}, the number $\gamma_j$ is equal to
  $$\dim_k V_{m-j-1}-\dim_k V_{m-j}.$$ Propositions \ref{prop-split} and \ref{prop-vanish} imply that the target of the morphism $d_i$ from
 \eqref{eq-connect} vanishes for every $i$ in $\{0,\ldots,m-1\}$, so that we can identify the $k$-vector space $V_i$ with
  $$H^0(\cC_k,\mathcal{L}_k\otimes\mathcal{F}_i).$$
    Using the description of the sheaves $\mathcal{F}_i$ in Propositions \ref{prop-split}, we can write
  $$\gamma_{j}=\dim_k H^0(\cC_{k,\red},\iota^*\mathcal{L}\otimes\mathcal{J}((j/m)\cC_{k})),$$ which is also equal to the Euler characteristic
  $$\chi(\cC_{k,\red},\iota^*\mathcal{L}\otimes\mathcal{J}((j/m)\cC_{k}))$$ by the vanishing result in Proposition \ref{prop-vanish}.
   Recall that the sheaf $\mathcal{J}(j/m)\cC_{k}$ was defined as the pullback of the line bundle $\mathcal{O}_{\cC}(\lfloor  (j/m)\cC_{k} \rfloor )$
    to the reduced strict normal crossings divisor $\sum_{i\in I_j}E_i$. Computing the above Euler characteristic as in \eqref{sss-resol}, we find that
   $$\gamma_j=(\sum_{i\in I_j}E_i\cdot \lfloor  (j/m)\cC_{k} \rfloor )+ \sum_{i\in I_j}g(E_i)-|I_j|+\sigma_j.$$

 It remains to compute the multiplicity of $0$ as a jump of $\Jac(C)$. By \eqref{sss-formula} and Lemma \ref{lemm-eldiv}, it is equal to
  the number of occurrences of the value $0$ in the tuple of elementary divisors of the inclusion of lattices
   $\Omega_{\log}(C)\subset \Omega_{\mathrm{can}}(C)$ (which is $g-u(C)$ by Proposition \ref{prop-term1})
plus the number  of occurrences of the value $m$ in the tuple of elementary divisors of the inclusion of lattices
 $$\Omega_{\log}(C)\otimes_R R^s \subset \Omega_{\mathrm{sat}}(C).$$ The latter number is equal to the dimension of the $k$-vector space
 $$V_{m-1}\cong H^0(\cC_k,\mathcal{L}_k\otimes\mathcal{F}_{m-1}),$$ but it follows from Proposition \ref{prop-split} that the sheaf $\mathcal{F}_{m-1}$ vanishes.
  Thus the multiplicity of $0$ as a jump of $\Jac(C)$ is equal to
  $g-u(C)$. By Proposition \ref{prop-split}, $g-u(C)$ is given by the expression in the statement for $j=0$ because $I_0=\{1,\ldots,r\}$.
\end{proof}

\begin{cor}
The jumps of $C$ only depend on the combinatorial reduction data of $C$, and not on the characteristic of $k$.
\end{cor}
\begin{proof}
This is obvious, since all the terms in the formula in Theorem \ref{theorem-formula} only depend on the combinatorial reduction data.
\end{proof}

\begin{cor}
The number of non-zero jumps of $\Jac(C)$ (counted with multiplicities) is equal to the unipotent rank of $\Jac(C)$.
\end{cor}
\begin{proof}
We have already remarked in the proof of Theorem \ref{theorem-formula} that the multiplicity of zero as a jump is $g-u(C)$, and the total number of jumps is $g$.
\end{proof}

\subsection{Jumps and principal components}

\sss Let $\cC$ be the minimal $sncd$-model of $C$ and write $\cC_k=\sum_{i=1}^rN_iE_i$ and $m=\lcm\{N_1,\ldots,N_r\}$ as before.  It is clear from Theorem \ref{theorem-formula}
 that all the jumps of $\Jac(C)$ are of the form $j/m=a/N_i$ for some $i$ in $\{1,\ldots,r\}$ and some $a$ in $\{0,\ldots,N_i-1\}$ (otherwise, the set $I_j$ is empty).
 However, we will now show that one can deduce much more
 precise information
 about which components $E_i$ can give rise to a jump. We keep the notations $I_j$ and $\sigma_j$ from Theorem \ref{theorem-formula}.

\begin{prop}\label{prop-lowbound}
 Let $j$ be an element of $\{1,\ldots,m-1\}$. If we denote by $\Gamma_j$ the dual graph of $\sum_{i\in I_j}E_i$ and by $\beta(\Gamma_j)$ its first Betti number, then the multiplicity of $j/m$ as a jump of $\Jac(C)$ is at least $$\beta(\Gamma_j)+\sum_{i\in I_j}g(E_i).$$
\end{prop}
\begin{proof}
We write
 $$(j/m)\cC_k=\lfloor  (j/m)\cC_{k} \rfloor+\{(j/m)\cC_{k}\}$$ as in the proof of Proposition \ref{prop-vanish}. The support of the $\Q$-divisor $\{(j/m)\cC_{k}\}$ is precisely the union of the components $E_i$ with $i \notin I_j$, and the coefficients in this $\Q$-divisor are strictly contained between $0$ and $1$.

 Let $D$ be a connected component of $\sum_{i\in I_j}E_i$. Then $ D \cdot \lfloor  (j/m)\cC_{k} \rfloor $ is strictly bigger than the negative of the number of intersection points of $D$ with $\sum_{i\notin I_j}E_i$. Summing over the connected components, we find that
 $$(\sum_{i\in I_j}E_i\cdot \lfloor  (j/m)\cC_{k} \rfloor )$$
 is at least $c_j-\sigma'_j$ where $c_j$ is the number of connected
 components of the divisor $\sum_{i\in I_j}E_i$ and $\sigma_j'$ the number of intersection points of $\sum_{i\in I_j}E_i$ with the remainder of the special fiber. Using the formula in Theorem \ref{theorem-formula} we find that the multiplicity of $j/m$ as a jump of $\Jac(C)$ is at least
 $$c_j-|I_j|+\sigma_j-\sigma'_j+\sum_{i\in I_j}g(E_i).$$ The value $\sigma_j-\sigma'_j$ is precisely the number of edges of the dual graph $\Gamma_j$ of $\sum_{i\in I_j}E_i$. Since
 $|I_j|$ is the number of vertices of $\Gamma_j$ and $c_j$ its number of connected components, we see that $c_j-|I_j|+\sigma_j-\sigma'_j$ equals the first Betti number of $\Gamma_j$.
\end{proof}

\begin{cor}\label{cor-posgen}
If $E$ is an irreducible component in $\cC_k$ of multiplicity $N$ and the genus of $E$ is at least one, then $a/N$ is a jump of $\Jac(C)$ for every $a$ in $\{1,\ldots,N-1\}$.
\end{cor}
\begin{proof} We set $j=am/N$ so that $j/m=a/N$.  By our assumption, there is at least one element $i\in I_j$ such that $g(E_i)>0$, so that $j/m$ is a jump of $\Jac(C)$ by Proposition \ref{prop-lowbound}.
\end{proof}

\begin{lemma}\label{lemma-integer}
Let $j$ be an element of $\{1,\ldots,N-1\}$ and let $i_0$ be an element of $I_j$. The intersection of $E_{i_0}$ with $\sum_{i\notin I_j}E_i$ is either empty
or consists of at least two points.
\end{lemma}
\begin{proof}
We once more use the equality
\begin{equation}\label{eq-integer}
(E_{i_0}\cdot \lfloor  (j/m)\cC_{k} \rfloor)=-(E_{i_0}\cdot \{(j/m)\cC_{k}\}).\end{equation}
Since the left hand side of \eqref{eq-integer} is an integer, we see that the intersection of $E_{i_0}$ with $\sum_{i\notin I_j}E_i$ cannot consist of a single point.
\end{proof}

\begin{prop}\label{prop-princ}
Each non-zero jump of $\Jac(C)$ is of the form $a/N$ where $N$ is the multiplicity of a principal component in $\cC_k$ and $a$ is an element of $\{1,\ldots,N-1\}$. Conversely, for every principal component of $\cC_k$ of multiplicity $N$, there exist a
 multiple $N'$ of $N$ and an element $a$ in $\{1,\ldots,N'-1\}$
 such that $a$ is prime to $N'$ and $a/N'$ is a jump of $\Jac(C)$.
\end{prop}
\begin{proof}
  Assume that $j$ is an element of $\{1,\ldots,m-1\}$ such that the components $E_i$ with $i\in I_j$ are all non-principal.
 We will show that $j/m$ cannot be a jump of $\Jac(C)$.
  A component $E_i$ with $i\in I_j$ cannot be a rational curve  intersecting the rest of $\cC_k$ in precisely one point, since otherwise, there would be a principal component
  $E_\ell$ with $\ell \in I_j$ by Lemma \ref{lemma-mult}. Thus each component $E_i$ with $i\in I_j$ is a rational curve intersecting the rest of $\cC_k$ in precisely two points. By Lemma \ref{lemma-integer} this is only possible if each connected component of $\sum_{i\in I_j}E_i$ consists of a single curve $E_i$ (note that $\sum_{i\in I_j}E_i$ cannot be
  the entire reduced special fiber $\cC_{k,\red}$ by our overall assumption that the multiplicities $N_i$ are coprime).
     By Theorem \ref{theorem-formula}, the multiplicity of $j/m$ as a jump of $\Jac(C)$ equals
  $$(\sum_{i\in I_j}E_i\cdot \lfloor  (j/m)\cC_{k} \rfloor ) -|I_j|+\sigma_j.$$
  Again using the equality \eqref{eq-integer}
   we see that the intersection product $$(\sum_{i\in I_j}E_i\cdot \lfloor  (j/m)\cC_{k} \rfloor )$$ is at most $-|I_j|$. On the other hand, $-|I_j|+\sigma_j$ is
   equal to $|I_j|$. Thus
      $$(\sum_{i\in I_j}E_i\cdot \lfloor  (j/m)\cC_{k} \rfloor ) -|I_j|+\sigma_j=0$$ and $j/m$ is not a jump of $C$.

Now we prove the converse statement.
 Let $i_0$ be an element of $\{1,\ldots,r\}$ such that $E_{i_0}$ is principal; we
 will prove that there exist a
 multiple $N'$ of $N_{i_0}$ and an element $a$ in $\{1,\ldots,N'-1\}$
 such that $a$ is prime to $N'$ and $a/N'$ is a jump of $\Jac(C)$.
 We may assume that there does not exist a principal component $E_{i_1}$ in $\cC_k$ such that $N_{i_1}$ is a multiple of $N_{i_0}$ and
 $N_{i_1}>N_{i_0}$ (otherwise we can simply replace $i_0$ by $i_1$). By Corollary \ref{cor-posgen} we may also assume that there does not exist
  a component of multiplicity $N_{i_0}$ with positive genus.
  Finally, if we set $j=m/N_{i_0}$, then by Proposition \ref{prop-lowbound}, we can also suppose that
  the divisor $\sum_{i\in I_j}E_i$ does not contain a loop. Otherwise, $\beta(\Gamma_j)$ is positive and $j/m=1/N_{i_0}$ is a jump.

   Assume that $\ell\in I_j$ and that $E_\ell$ is a non-principal component of $\cC_k$.
    Then $E_\ell$ must intersect the rest of $\cC_k$ in exactly two points, because otherwise there would exist a principal component in $\cC_k$ whose multiplicity is a strict multiple of
  $N_{i_0}$ by Lemma \ref{lemma-mult}, contradicting the maximality of $N_{i_0}$.
  If $E$ and $E'$ are the components of $\cC_k$ intersecting $E_\ell$  (where possibly $E=E'$) then Lemma \ref{lemma-integer} implies that either
  $E$ and $E'$ are both contained in $\sum_{i\in I_j}E_i$, or neither of them is. From these observations we deduce that
 some connected component of $\sum_{i\in I_j}E_i$ must intersect $\sum_{i\notin I_j}E_i$ in at least three points.
  We denote by $M_1$, $M_2$ and $M_3$ the multiplicities
   of the components of $\sum_{i\notin I_j}E_i$ intersecting $\sum_{i\in I_j}E_i$ for an arbitrary choice of three such intersection points.
     Then for at least one of both elements $a$ in $\{1,N_{i_0}-1\}$ we have
     $$\{aM_1/N_{i_0}\}+\{aM_2/N_{i_0}\}+\{aM_3/N_{i_0}\} <2,$$ where $\{q\}=q-\lfloor q\rfloor$ denotes the fractional part of a rational number $q$.
     In the notation of the proof of Proposition \ref{prop-lowbound},
   this implies that
   $$(\sum_{i\in I_{aj}}E_i\cdot \lfloor  (aj/m)\cC_{k} \rfloor )$$ is at least $c_j+1-\sigma'_j$, and the arguments in that proof show that
   $aj/m=a/N_{i_0}$ is a jump of $\Jac(C)$.
\end{proof}

\begin{cor}\label{cor-lcm}
The stabilization index $e(C)$ of $C$ is equal to least common multiple of the
 denominators of the jumps of $\Jac(C)$, that is, the smallest positive integer $e$ such that $e\cdot j$ is an integer for every
 jump $j$ of $\Jac(C)$. Thus the stabilization index of $C$ is equal to the stabilization index of $\Jac(C)$ in the sense of \eqref{sss-stabdef}.
\end{cor}
\begin{proof}
This is an immediate consequence of Proposition \ref{prop-princ}.
\end{proof}

\begin{cor}\label{cor-stab} The stabilization index $e(C)$
 is equal to the smallest possible saturation index of an $R$-model $\cC$ of $C$ such that $\cC^+$ is log regular.
\end{cor}
\begin{proof}
 If $m$ is the saturation index of such a model $\cC$, then it follows from Theorems \ref{thm-logtame} and \ref{thm-logdescr} that
  $m\cdot j$ is an integer for every jump $j$ of $C$ (this was already observed in the proof of Theorem \ref{thm-rat}).
  Thus by Corollary \ref{cor-lcm}, $m$ must be divisible by $e(C)$. On the other hand, we have constructed in the proof of Theorem \ref{thm-rat}
  a model $\cC$ such that $\cC^+$ is log regular and such that the saturation index of $\cC^+\to S^+$ is precisely $e(C)$.
\end{proof}

\subsection{Examples}
\sss We illustrate with a few examples how the formula in Theorem \ref{theorem-formula} can be used to compute the jumps in practice.
 We start with the case where $C$ is an elliptic curve. Then $C$ has precisely one jump, which is zero if and only if $C$ has semi-stable reduction. If this is not the case, the unipotent rank $u(C)$ equals one, and the unique jump is non-zero.

 Let $\cC$ denote the minimal $sncd$-model of $C$. If $C$ has reduction type $I^*_{\geq 0}$, the multiplicity of any irreducible component of $\cC_k$ is either $1$ or $2$, and the jump is therefore $1/2$. In all other remaining cases, the special fiber $\cC_k$ has precisely one principal component, which is a rational curve
  intersecting the other components in precisely three points. We denote its multiplicity by $N$.
 Then we know by Proposition \ref{prop-princ} that the jump of $C$ is of the form $a/N$. The proof of Proposition \ref{prop-princ} even tells us how to find
 $a$: it is the unique element in $\{1,N-1\}$ that satisfies the property that
 $$\{aM_1/N\}+\{aM_2/N\}+\{aM_3/N\} <2$$ where $M_1$, $M_2$ and $M_3$ are the multiplicities of the components of $\cC_k$ intersecting the principal component.
  In this way, we immediately recover the table from \cite[5.4.5]{edix} and \cite[Table~8.1]{Halle-neron} giving the jump for each of the Kodaira-N\'eron reduction types (see Table \ref{table-ell}).

 \begin{table}[h!]\label{table-ell}
\begin{tabular*}{0.8\textwidth}
     {@{\extracolsep{\fill}}
     |c||c|c|c|c|c|c|c|c|} \hline     & & & & & & & & \\  Type & $I_{\geq 0}$ &
$II$ & $III$ & $IV$ & $I^*_{\geq 0}$ & $IV^*$ & $III^*$ & $II^*$
\\[1ex] \hline & & & & & & & &
\\ jump & 0  &1/6 & 1/4 & 1/3 & 1/2 & 2/3 & 3/4 & 5/6
\\[1ex] \hline
\end{tabular*}
\\[2ex]
\refstepcounter{subsubsection} \caption{The jump of an elliptic curve} \label{table-elliptic}
\end{table}

 \sss Now we consider the case where $C$ has genus 2. Let $\cC$ be the minimal $sncd$-model of $C$ and assume
 that $\cC_k$ has a component of genus one and multiplicity $N>1$. Then the unipotent rank of $\Jac(C)$ is one, so that $0$ is a jump of $C$ of multiplicity one.
  On the other hand, Corollary \ref{cor-posgen} tells us that $a/N$ is a jump of $C$ for every $a$ in $\{1,\ldots,N-1\}$. Since the total number of jumps of $C$ is equal to $2$,
   we must have $N=2$, and the jumps of $C$ are $0$ and $1/2$.
    With similar arguments one can quickly reproduce the tables in \cite[\S8.4]{Halle-neron} which give the
   jumps of $C$ for each of the  reduction types in the Namikawa-Ueno classification of degenerations of genus 2 curves.

\bibliographystyle{numsty}

\end{document}